\input amstex
\documentstyle{amsppt}

\define\ZZ{\bold{Z}}
\define\QQ{\bold{Q}}
\define\RR{\bold{R}}
\define\CC{\bold{C}}
\define\ee{\bold{e}}
\define\oo{\infty}
\define\a{\alpha}
\define\be{\beta}
\define\ep{\epsilon}
\define\Gm{\Gamma}
\define\de{\delta}
\define\pr{\prime}
\define\s{\sigma}
\define\ar{\text {arg}}
\define\ti{\times}
\define\lr{\longrightarrow}
\define\la{\langle}
\define\ra{\rangle}
\define\vt{\vartheta}
\define\Ga{\text {Gal}}
\define\sumprime{\sideset \and^{\prime} \to\sum}

\NoRunningHeads

\magnification=1200
\hoffset=0.25in

\topmatter
\title
On calculations of zeros of various $L$-functions
\endtitle
\author
By Hiroyuki Yoshida
\endauthor
\endtopmatter

\document

\footnote""{${}^*$ During the final stage of writing this paper, the author was at
MSRI supported in part by NSF grant $\sharp$DMS9022140.}

Introduction.  As we have shown several years ago [Y2], zeros of $L(s, \Delta )$
and $L^{(2)}(s, \Delta )$ can be calculated quite efficiently by a certain experimental
method. Here $\Delta$ denotes the cusp form of weight $12$ with respect to
$SL(2, \bold Z)$ and $L(s, \Delta )$ (resp. $L^{(2)}(s, \Delta )$) denotes
the standard (resp. symmetric square) $L$-function attached to $\Delta$. 
The purpose of this paper is to show that this method can be applied to 
a wide class of $L$-functions so that we can obtain precise numerical values of their 
zeros. \footnote {After the publication of [Y2], H. Ishii [Is] published a table
of zeros of standard $L$-functions attached to modular forms for $15$ cases.
It also comes to the author's notice that a program of the calculation of
zeros of $L(s, \Delta )$ is included in \lq\lq Mathematica \hskip -0.3em \rq\rq \, package, following
the method of [Y2].}

We organize this paper as follows. In \S1, we shall describe basic features of our method
of calculation, which is repeated applications of partial summation. In \S2,
we shall study the $r$-th symmetric power $L$-function $L^{(r)}(s, \Delta )$
attached to $\Delta$. Since the cases $r=1$, $2$ are discussed in [Y2], we shall
exclusively treat the cases $r=3$, $4$. In \S3, we shall study the $L$-functions
attached to modular forms of half integral weight. These $L$-functions do not have 
Euler products. Naturally the Riemann hypothesis fails for them; we shall find 
many zeros off the critical line, though major part of zeros lie on the critical line. 
We shall also calculate the location of these zeros off the critical line. Though there is some
hope to find relations among zeros of $L$-functions of two modular forms which are in
the Shimura correspondence, no explicit results came out so far.

In \S4, we shall study $L$-functions attached to Hecke characters of non-$A_0$ type 
of real quadratic fields. D.A. Hejhal showed great interest to make experiments in this case, 
since coefficients are non-computable combinatorially;
hence there is a slight possibility that the Riemann hypothesis may break down
for these $L$-functions.
We have made experiments on $44$ cases summarized in Table 4.3; so far
no counterexamples are found.

In \S5, we shall study the Artin $L$-function attached to a $4$-dimensional
non-monomial representation of $\text {Gal} \, (\overline {\bold Q}/\bold Q)$.
In \S6, we shall discuss the controle of error estimates in our calculation. In \S7, we shall
consider the explicit formula for the $L$-function attached to a modular form
of weight $8$ with respect to $\Gamma _0(2)$. We shall compare both sides
of the explicit formula numerically. In \S8, we shall present sample programs to 
compute values of $L$-functions, which may be convenient for the reader. In \S9,
we shall formulate a conjecture which emerged during the process of our experiments.

Most of sections have attached tables to show results explicitly. Concerning actual 
computations, we have used \lq\lq UBASIC \hskip -0.2em\rq\rq \, created by Y. Kida.
(It was not available when we wrote [Y2].) The calculation was done
by personal computers which are not necessarily so fast. However our experiments 
extended over long time (about three years) and UBASIC is quite fast
(compared with some other softwares) for numerical calculations, the author thinks
that our tables are fairly extensive.

A motivation in these calculations has been to find non-trivial \lq\lq functorial \hskip -0.3em \rq\rq
\, prope-\linebreak rties which may exist among zeros of $L$-functions, as was hinted
in [Y1]. Though our experiments are not successful in this regard, conjectures
stemmed from them are formulated in \S9.

We can pursue these calculations still further. The topics which may be
included in this paper are:
\item {1)} The Hasse-Weil zeta functions of algebraic curves, for example
$y^2=x^5-x+1$.
\item {2)} The Dirichlet series
$\dsize \sum _{n=1}^{\infty} \frac {n\alpha -\lbrack n\alpha \rbrack -1/2}{n^s}$
studied by Hecke [H], where $\alpha$ is a real irrational number.
\item {3)} Applications to Riemann-Siegel type formulas.
\item {4)} Calculations of critical values of $L$-functions.

\noindent Our results on these topics are still fragmentary, so the full discussion should
be postponed to future occasions.

\bigskip

Notation. For a complex number $z$, we denote by $\Re(z)$ (resp. $\Im(z)$) the real
(resp. imaginary) part of $z$. The letter $q$ stands for $\exp (2\pi \sqrt {-1}z)$
when it is clear from the context. For modular forms, we follow the notation in
Shimura [Sh1].

\bigskip

\centerline {\S1. An overview on our method of calculations}

\medskip

Let
$$
L(s)=\sum _{n=1}^{\oo} a_nn^{-s}
\tag{1.1}
$$
be a Dirichlet series which is absolutely convergent when
$\Re (s)>\sigma$ for some $\sigma >1$. In this paper, we shall consider only such
$L(s)$ which can be analytically continued to the whole complex plane as an entire function 
and satisfies a functional equation of the form
$$
R(k- \bar s)=\kappa \overline {R(s)} .
\tag{1.2}
$$
Here $\kappa$ is a constant of absolute value $1$, $k>0$,
$$
R(s)=N^s \prod _{i=1}^m \Gm (b_is+c_i)L(s)
$$
with $N>0$, $b_i>0$, $c_i \in \CC$. We note that (1.2) is equivalent to
$$
R(k- s)=\kappa R(s) 
\tag{1.3}
$$
if $a_n \in \RR$, $c_i \in \RR$ for all $n$ and $i$. Put
$\kappa ^{-1}=\kappa _1^2$ with some $\kappa _1 \in \CC$. By (1.2), we have
$$
\kappa _1R(s) \in \RR \qquad \text {if} \quad \Re (s)=k/2.
\tag{1.4}
$$

Take any $\de >0$. For $T>0$, let $N(\de ,T)$ denote the number
of zeros of $R(s)$ counted with multiplicity in the domain
$$
-\de \leq \Re (s) \leq k+\de , \qquad 0 \leq \Im (s) \leq T .
$$
Let $D$ be the rectangle whose vertices are $-\de$, $k+\de$, $k+\de +iT$,
$-\de +iT$ and let $C$ denote the contour $\partial D$ taken in positive direction.
By the argument principle, we have
$$
N(\de ,T)= \frac {1}{2\pi i} \int _C \frac {R^{\pr}(s)}{R(s)} \, ds ,
\tag{1.5}
$$
assuming that neither zeros nor poles of $R(s)$ lie on $C$. Let $C_1$ denote
the portion of $C$  from $k/2$ to $k/2+iT$. By the functional equation (1.2),
we obtain
$$
N(\de ,T)= \pi ^{-1} \Delta \ar \ R(s)
=\pi ^{-1} \Delta (\ar \ N^s \prod _{i=1}^m \Gm (b_is+c_i)) +
\pi ^{-1} \Delta (\ar \ L(s)) ,
\tag{1.6}
$$
where $\Delta \ar$ denotes the variation of the argument on $C_1$, i.e.,
from $s=k/2$ to $k/2+iT$ along $k/2$ to $k+\de$, $k+\de$ to $k+\de +iT$,
$k+\de +iT$ to $k/2+iT$. Set
\footnote {When it is clear from the context, we shall use $\vt (T)$ for
the \lq\lq phase factor" of this type in the following sections without
further explanation.}
$$
\vt (T)= \Delta \ar \ (N^s \prod _{i=1}^m \Gm (b_is+c_i))  .
$$
Assume 
$$
\Re (c_i) >-b_ik/2 \qquad \text {for} \quad 1 \leq i \leq m.
$$
Then since $b_i>0$, $N^s \prod _{i=1}^m \Gm (b_is+c_i)$
has neither zeros nor poles in the domain $\Re (s) \geq k/2$. Hence $\vt (T)$
is equal to the variation of the argument of $N^s \prod _{i=1}^m \Gm (b_is+c_i))$
on the line segment $[ k/2, k/2+iT]$. We note that $\vt (T)$ can be
computed in high precision very easily using Stirling's formula (cf. [WW], p. 252)
combined with the relation $\Gm (s+1)=s\Gm (s)$. We obtain
$$
N(\de ,T)=\pi ^{-1} \vt (T)+\pi ^{-1} \Delta (\ar \ L(s)).
\tag{1.7}
$$

Now let us consider the case when $R(s)$ has zeros in $(-\de , k+\de )$.
Let $r$ denote the number of zeros of $R(s)$ i.e., of $L(s)$, in this interval
counted with multiplicity. Then (1.7) holds with the modification
$$
N(\de ,T) -\frac {r}{2} =\pi ^{-1} \vt (T)+\pi ^{-1} \Delta (\ar \ L(s)).
\tag{1.8}
$$
Here $\Delta (\ar \ L(s))$ is counted by dividing $C_1$ into a finite number of
paths removing real zeros of $L(s)$ and summing the variations of the argument
of $L(s)$ on each of them. The validity of (1.8) can be seen by modifying
$C$ by small semi-circles which detour the real zeros of $L(s)$.

Throughout the paper, to compute $L(s)$, we shall employ our method
given in [Y2], which is repeated applications of Abel's partial summation.
Set
$$
s_n^{(0)} = a_n, \qquad  u_n^{(0)} = n^{-s}
$$
and define $s_n^{(l)}$, $u_n^{(l)}$ recursively by
$$
s_n^{(l)} = \sum\limits_{m=1}^n s_m^{(l-1)}, \qquad
u_n^{(l)} = u_n^{(l-1)} - u_{n+1}^{(l-1)}, \qquad l \ge 1 .
\tag{1.9}
$$
Put $S_N^{(l)} = \sum\limits_{n=1}^N s_n^{(l)}u_n^{(l)}$. Then we have
$$
S_N^{(l)} = S_N^{(l-1)} - s_N^{(l)}u_{N+1}^{(l-1)}.
\tag{1.10}
$$
As we have seen in [Y2], in several cases, $S_N^{(l)}$ seems to approximate
$L(s)$ amazingly well when we choose $N$ and $l$ sufficienly large. In the
succeeding sections, we shall present various types of $L$-functions which
can be treated in more or less similar fashion. The efficacy of our method
seems to depend strongly on the arithmetical nature of the coefficients $a_n$
of a Dirichlet series $L(s)$.

We shall conclude this section by technical remarks concerning actual
computations of $S_N^{(l)}$. As the first step, we should construct a table
of $a_n$. For Dirichlet series considered in this paper, this step can be
achieved rather easily. Since we can compute $S_N^{(l)}$ from $S_N^{(0)}$
by (1.10), the computation of $S_N^{(0)} = \sum\limits_{n=1}^N a_n n^{-s}$
is the substantial and the most time consuming part of our calculation. However
usually $s_n^{(l)}$ becomes very large and $u_n^{(l)}$ very small when $l$
increases. Therefore it is indispensable to perform the actual computation in
high precision. For $u_n^{(l)}$, the following formula (1.11) should preferably
be used than to compute it directly from the definition.
$$
\aligned
u_N^{(l)}&=N^{-s} \sum _{k=l}^{\oo}
(\sum _{m=1}^l (-1)^m \binom lm m^k) \\
&(-1)^k \frac {s(s+1) \cdots (s+k-1)} {k !} N^{-k}, 
\endaligned
\qquad \qquad N>l \ge 1.
\tag{1.11}
$$
If we replace $\sum _{k=l}^{\oo}$ by $\sum _{k=l}^L$, the error is less than
$2^l \frac { \vert s (s+1) \cdots (s+L) \vert } {(L+1)!} 
(\frac {l}{N} )^{L+1}
\vert N^{-s} \vert$ if $\Re (s) \ge -L-1$.

\newpage

\centerline
{\S2. $L$-functions attached to symmetric tensor representations of $GL(2)$}

\medskip

Let $f(z)= \sum _{n=1}^{\oo} c_ne^{2\pi inz} \in S_k(SL(2, \ZZ ))$ be
a normalized common eigenfunction of Hecke operators. The $L$-function
$L(s,f)=\sum _{n=1}^{\oo} c_n n^{-s}$ attached to $f$ converges absolutely
when $\Re (s)> \frac {k+1}{2}$ and has the Euler product
$$
L(s,f)= \prod _p (1-c_pp^{-s} +p^{k-1-2s})^{-1} .
$$
Put
$$
1-c_pX+p^{k-1}X^2=(1-\a _pX)(1-\be _pX)
$$
with $\a _p$, $\be _p \in \CC$, where $X$ is an indeterminate. For a positive
integer $r$, we define an Euler product
$$
L^{(r)}(s, f)= \prod _p [ (1- \a _p^r p^{-s})(1-\a _p^{r-1}\be _p p^{-s})
\cdots (1-\be _p^r p^{-s}) ]^{-1}
\tag{2.1}
$$
which converges absolutely when $\Re (s)> \frac {r(k-1)}{2} +1$.
It is conjectured that $L^{(r)}(s, f)$ can be analytically continued
to the whole complex plane as an entire function and satisfies a functional
equation. The conjectural functional equation of $L^{(r)}(s,f)$ takes the
following form (cf. Serre [Se]). If $r$ is odd, put $r=2m-1$,
$$
R^{(r)}(s,f)=(2\pi )^{-ms} \prod _{i=0}^{m-1}
\Gm (s-i(k-1))L^{(r)}(s,f),
\tag{2.2}
$$
$$
\ep _r=(\sqrt {-1})^{m+(k-1)m^2}.
\tag{2.3}
$$
If $r$ is even, put $r=2m$,
$$
\aligned
R^{(r)}(s,f)= \pi ^{-s/2}(2\pi )^{-ms}& (\prod _{i=0}^{m-1}
\Gm (s-i(k-1))) \\
&\Gm (\frac {s-m(k-1)+\delta }{2}) L^{(r)}(s,f),
\endaligned
\tag{2.4}
$$
$$
\ep _r=(\sqrt {-1})^{m+(k-1)m(m+1)+\delta },
\tag{2.5}
$$
where $\delta =0$ (resp. $1$) if $m$ is even (resp. odd). Then the
functional equation
$$
R^{(r)}(s,f)=\ep _r R^{(r)}(r(k-1)+1-s, f)
\tag{2.6}
$$
is predicted. 
A quick way to see (2.6) is as follows. Let $M_f$ be the motive of rank $2$
over $\QQ$ attached to $f$. We see that the Hodge realization of $M_f$
corresponds to the two dimensional representation
$$
\rho = \text {Ind} (\psi _k ; W_{\CC} \lr W_{\RR , \CC } )
$$
of $W_{\RR , \CC }$. Here $W_{\CC}= \CC ^{\times}$, $W_{\RR, \CC}$
is the Weil group of $\CC$ over $\RR$ and $\psi _k$ is the quasi-character
$\psi _k(x)=x^{-(k-1)}$ of $W_{\CC}$. Let 
$\s _r : GL(2) \lr GL(r+1)$ be the symmetric tensor representation of degree
$r$ and put $\rho _r=\s _r \circ \rho$. Then we find
$$
\rho _r \cong \oplus _{i=0}^{m-1} \text {Ind}
(x \lr x^{-(r-i)(k-1)} \bar x ^{-i(k-1)} ; \ W_{\CC} \lr W_{\RR, \CC}),
\quad r=2m-1,
\tag{2.7}
$$
$$
\aligned
\rho _r \cong &\oplus _{i=0}^{m-1} \text {Ind}
(x \lr x^{-(r-i)(k-1)} \bar x ^{-i(k-1)} ; \ W_{\CC} \lr W_{\RR, \CC}) \\
&\oplus \{ (x \lr \vert x \vert ^{-m(k-1)}( \text {sgn} \ x)^{m(k-1)}) \circ t \},
\quad r=2m,
\endaligned
\tag{2.8}
$$
where $t$ denotes the transfer map from $W_{\RR , \CC}$ to
$\RR ^{\times}$. The gamma factor and the constant $\ep _r$ of the functional
equation can be calculated as the usual gamma factor and the constant attached to
the representation $\rho _r$ of $W_{\RR , \CC}$; hence we obtain
(2.2) $\sim$ (2.6).

We refer the reader to Shahidi [Sha1], [Sha2] for what are known 
on these symmetric power $L$-functions, in more general cases.

Let $\Delta (z)=q \prod _{n=1}^{\oo} (1-q^n)^{24} \in S_{12}(SL(2, \ZZ ))$,
$q=e^{2\pi \sqrt {-1} z}$. The calculation of zeros of $L^{(r)}(s, \Delta )$
for $r=1$, $2$ is given in [Y2]. We consider the case $r \ge 3$. To compute
$L^{(r)}(s, \Delta )$, we modify our summation method slightly in the
following way. Fix $r$, choose $v=v_r>0$ and set
$$
L^{(r)}(s, \Delta )=\sum _{n=1}^{\oo} a_n n^{-s}
=\sum _{n=1}^{\oo} (a_n n^{-v}) n^{-(s-v)} .
$$
Put
$$
s_n^{(0)}=a_n n^{-v}, \qquad u_n^{(0)}=n^{-(s-v)}
$$
and define $s_n^{(l)}$, $u_n^{(l)}$ recursively by (1.9). We set
$S_N^{(l)}= \sum _{n=1}^N s_n^{(l)} u_n^{(l)}$. It turns out that
a suitable choice of $v$ depending on $r$ yields good results.
We can interpret this as the neutralization of the effect of extremely
large value of $s_n^{(l)}$ and extremely small value of $u_n^{(l)}$.

As the first example, let $r=3$. We take $v=8$. For $s=17+it$, $t=20$,
the values of
$$
R_j=\Re (\exp (i\vt (t)) S_N^{(j)}), \qquad
I_j=\Im (\exp (i\vt (t)) S_N^{(j)})
$$
are given in Table 2.1. In Table 2.2, we give the values of $t_n$ the $n$-th
zero of $L^{(3)}(s, \Delta )$, $s=17+it$ on the critical line for $0 \le t \le 40$.

Next we apply our summation method to $L^{(4)}(s, \Delta )$ taking $v=12$.
For $s=\frac {45}{2}+it$, $t=10$, the values of
$$
R_j=\Re (\exp (i\vt (t)) S_N^{(j)}), \qquad
I_j=\Im (\exp (i\vt (t)) S_N^{(j)})
$$
are given in Table 2.3. In Table 2.4, we give the values of $u_n$ the $n$-th
zero of $L^{(4)}(s, \Delta )$, $s=\frac {45}{2}+iu$ on the critical line for 
$0 \le u \le 20$.

We can see, by the same technique as will be given in \S3 and \S4, that
the Riemann hypothesis holds for $L^{(3)}(s, \Delta )$ (resp. $L^{(4)}(s, \Delta )$)
in the range $0 \leq \Im (s) \leq 40$ (resp. $0 \leq \Im (s) \leq 20$ ) and that
the zeros $17+it_n$ (resp. $\frac {45}{2} +iu_n$) are simple.

\bigskip

\centerline
{\S3. Modular forms of half integral weight}

\medskip

Put
$$
\theta (z)=\sum _{n \in \ZZ} \ee (n^2z)=1+2 \sum _{n=1}^{\oo} q^{n^2}, 
\qquad \eta (z)=\ee (z/24) \prod _{n=1}^{\oo} (1- \ee (nz)) .
$$
By Shimura, [Sh2], (4.1), we have
$$
\dim S_8(\Gm _0(2))=1, \qquad \dim S_{9/2}(\Gm _0(4))=1
$$
and $(\eta (z)\eta (2z))^8$ (resp. $\theta (z)^{-3}\eta (2z)^{12}$) spans
$S_8(\Gm _0(2))$ (resp. $S_{9/2}(\Gm _0(4))$). Put
$$
\gather
f(z)=(\eta(z) \eta (2z))^8=\sum _{n=1}^{\oo} a_nq^n, \qquad
g(z)=\theta (z)^{-3}\eta (2z)^{12}=\sum _{n=1}^{\oo} c_nq^n, \\
L(s,f)=\sum _{n=1}^{\oo} a_nn^{-s}, 
\qquad L(s,g)=\sum _{n=1}^{\oo} c_nn^{-s}, \\
R(s,f)=2^{s/2}(2\pi )^{-s}\Gm (s)L(s,f), 
\qquad R(s,g)=2^s(2\pi )^{-s}\Gm (s)L(s,g).
\endgather
$$
Then $f$ and $g$ are in the Shimura correspondence; $L(s,f)$ and $L(s,g)$
can be analytically continued to the whole complex plane as entire functions
and satisfy the functional equations
$$
R(s,f)=R(8-s,f), \qquad R(s,g)=R(9/2-s, g) .
\tag{3.1}
$$
This example is described in detail in [Sh2], \S4.
For $t>0$, let $\vt _f(t)$ (resp. $\vt _g(t)$) denote the variation of the argument
of $2^{s/2}(2\pi )^{-s} \Gm (s)$ (resp. $2^s(2\pi )^{-s} \Gm (s)$) from $4$ to
$4+it$ (resp. $9/4$ to $9/4+it$). 

For $L(s,f)$, $s=4+it$, $t=100$, the values of
$$
R_j=\Re (\exp (i\vt _f(t)) S_N^{(j)}), \qquad
I_j=\Im (\exp (i\vt _f(t)) S_N^{(j)})
$$
are given in Table 3.1.

For $L(s,g)$, $s=\frac {9}{4}+it$, $t=100$, the values of
$$
R_j=\Re (\exp (i\vt _g(t)) S_N^{(j)}), \qquad
I_j=\Im (\exp (i\vt _g(t)) S_N^{(j)})
$$
are given in Table 3.2.

By our method, we can compute zeros of $L(s,f)$ and of $L(s,g)$ on the critical line
with sufficient accuracy observing sign changes of $e^{i\vt _f(t)}L(4+it,f)$ and
$e^{i\vt _g(t)}L(\frac {9}{4}+it,g)$. In Table 3.3, we list the $n$-th zero $t_n$
of $L(s,f)$, $s=4+it$ in the range $0 \le t \le 100$. In Table 3.4, we list the $n$-th zero
$u_n$ of $L(s,g)$, $s=\frac {9}{4}+iu$ in the range $0 \le u \le 100$.

Now let us examine the Riemann hypothesis for $L(s,f)$. We see $f(iy)>0$ for
$y>0$ by the product expansion of the $\eta$-function. By the integral
representation
$$
(2\pi )^{-s}\Gm (s)L(s,f)= \int _0^{\oo} f(iy)y^{s-1} \, dy ,
$$
we see that $L(s,f)>0$ for $s>0$. For $T>0$, let $N(T)$ denote the number of
zeros of $L(s,f)$ counted with multiplicity in the domain
$$
\vert \Re (s) -4 \vert <1/2, \qquad 0 \le \Im (s) \le T.
$$
By (1.7) taking $\de =1/2$, we have
$$
N(T)=\pi ^{-1}\vt _f(T)+\pi ^{-1} \Delta\text {arg} (L(s,f)).
\tag{3.2}
$$
Since $L(s,f) \ne 0$ if $\Re (s) -4 \ge 1/2$, $\Delta\text {arg} (L(s,f))$
equals the variation of the argument of $L(s,f)$ along the line segments
$L_1=[4, 4+\frac{1}{2}+\mu ]$, $L_2=[4+\frac{1}{2}+\mu , 4+\frac{1}{2}+\mu +iT]$,
$L_3=[4+\frac{1}{2}+\mu +iT, 4+iT]$ for any $\mu >0$. Take $T=100$, $\mu =1$.
Then we have $\pi ^{-1}\vt _f(100)=69.0171 \cdots$. Hence if we can show
$\vert \Delta \text {arg} (L(s, f)) \vert < \pi /2$, we can conclude $N(T)=69$.
For this purpose, it suffices to show that $\Re (L(s,f))>0$ when $s \in L_i$,
$i=1$, $2$, $3$. For $L_1$, we have shown this fact above. For $L_2$, this
fact can be proved as in [Y2], \S4. For $L_3$, we divide it into $150$ small
intervals and appeal to our heuristic calculation. We have observed
$$
\Re (L(s,f))>0.83 \qquad \text {on} \quad L_3.
$$
Thus we conclude that the Riemann hypothesis holds for $L(s,f)$ when
$0 \le \Im (s) \le 100$. All the zeros are simple.

Now let us consider zeros of $L(s,g)$. We have $g(iy)>0$ for $y>0$
since $\theta (iy)>0$, $\eta (iy)>0$ for $y>0$. By the integral representation
$$
(2\pi )^{-s}\Gm (s)L(s,g)= \int _0^{\oo} g(iy)y^{s-1} \, dy ,
$$
we see that $L(s,g)>0$ for $s>0$. For $T>0$, let $N(T)$ denote
the number of zeros of $L(s,g)$ counted with multiplicity in the domain
\footnote {The zero free region of $L(s,g)$ is non-trivial. Here we content
ourselves by regarding $\vert \Re (s) -9/4 \vert \leq 2$ is \lq\lq sufficiently wide".}
$$
\vert \Re (s) -9/4 \vert \le 2, \qquad 0 \le \Im (s) \le T.
$$
By (1.7) taking $\de =2$, we have
$$
N(T)=\pi ^{-1}\vt _g(T)+\pi ^{-1} \Delta\text {arg} (L(s,g)).
$$
Here $\Delta\text {arg} (L(s,g))$ denotes the variation of the argument of $L(s,g)$ 
along the line segments $L_1=[9/4, 17/4]$, $L_2=[ 17/4 , 17/4+iT]$,
$L_3=[ 17/4+iT, 9/4+iT]$. Take $T=100$. Then we have 
$\pi ^{-1}\vt _g(100)=79.1885 \cdots$. Dividing $L_2$ and $L_3$ into small
intervals, we have observed $\pi ^{-1} \Delta (\text {arg} L(s,g))= 0.8114 \cdots$. 
Thus we obtain $N(T)=  80$. On the otherhand, we have obtained only $54$ zeros
on the critical line. Therefore, assuming that these zeros are simple, there must 
exist $13$ zeros in the right-hand side of the critical line:
$9/4 < \Re (s) \le 17/4$, $0 \le \Im (s) \le 100$. These zeros, together with those in
$100 \le \Im(s) \le 150$, are given in Table 3.5.

Our method of calculation of these exceptional zeros is as follows. Let us
consider a box $B$ given by $d_1 \le \Re (s) \le d_2$, $h_1 \le \Im (s) \le h_2$.
By the argument principle, we can determine whether $L(s,g)$ has a zero
incide $B$ or not. First we find a box $B$ in which $L(s,g)$ has zeros
by trial and error. Then dividing $B$ into sub-boxes and applying the
principle above successively, we can obtain a good approximation for a zero
inside of $B$.\footnote {It is more efficient to use a variant of the Newton
method once we get a rough approximation. After finding a precise location,
final check should be done by the method described above.}

The Riemann hypothesis does not hold for $L(s,g)$. This should be of
no surprise since $L(s,g)$ does not have an Euler product.

We shall study one more example of modular forms of half integral weight.
For $k \ge 1$, put
$$
G_k(z)=\frac {1}{2} \zeta (1-k) + \sum _{n=1}^{\oo} \s _{k-1} (n)q^n,
\qquad \s _{k-1}(n)= \sum _{d \mid n, d>0} d^{k-1} .
$$
Let
$$
\de (z)= \frac {60}{2\pi i} (2G_4(4z)\theta ^{\pr} (z) -G_4^{\pr} (4z)\theta (z))
=\sum _{n=1, n \equiv 0,1 \mod 4}^{\oo} c(n)q^n .
$$
Then we have $\de \in S_{13/2}(\Gm _0(4))$ (cf. Kohnen-Zagier [KZ], p. 177)
and $\de$ corresponds to $\Delta (z)=\eta (z)^{24} \in S_{12}(SL(2, \ZZ ))$
under the Shimura correspondence. The values of $c(n)$ can easily be computed by
$$
c(n)=\omega (\sqrt {n}) \cdot n + 120 \sum _{m=1}^{[\frac {n-1}{4}]}
\omega (\sqrt {n-4m}) \s _3(m) (2n-9m) -15n \s _3(n/4),
$$
where
$$
\omega (x)= \cases
1 \qquad \text {if} \quad x \in \ZZ , \\
0 \qquad \text {if} \quad x \notin \ZZ .
\endcases
$$
Let
$$
\de _0(z)= \sum _{n \equiv 0 \mod 4} c(n)q^{n/4} .
$$
Then we have $\de _0 \in S_{13/2} (\Gm _0(4))$ and
$$
\de (\frac {i}{4y}) = -\sqrt {2} y^{6+1/2} \de _0 (iy), \qquad y>0.
$$
(cf. [KZ], p. 190.) Put
$$
h_{\pm} (z)=(\de \mp 2^{-6}\de _0)(z), \qquad
R(s, h_{\pm})=2^s (2\pi )^{-s} \Gm (s)L(s, h_{\pm}).
$$
Then we have
$$
h_{\pm}( \frac {i}{4y}) = \pm (2y)^{6+1/2} h_{\pm}(iy), \qquad y>0.
$$
Hence we obtain the functional equations for the entire functions $R(s, h_{\pm})$:
$$
R(\frac {13}{2} -s, h_{\pm}) =\pm R(s, h_{\pm}).
$$

We have computed zeros of $L(s, h_{\pm})$ on the critical line $\Re (s)=13/4$
in the range $0 \le \Im (s) \le 100$. The results are given in Table 3.6. There 
$u_n^{\pm}$ denotes the $n$-th zero of $L(s, h_{\pm})$ for $s=\frac {13}{4} +iu$.

In Table 3.7, we also list zeros of $L(s, h_{\pm})$ not on the critical line
which are searched in the region $13/4 < \Re (s) \le 25/4$, $0 \le \Im (s) \le 100$;
$\rho _i^{\pm}$ denotes a zero of $L(s, h_{\pm})$. It is a very interesting phenomenon
that $L(s,h_{-})$ has much more zeros outside of the critical line compared with
$L(s,h_{+})$.

\bigskip

\centerline
{\S4. $L$-functions attached to Hecke characters of infinite order of real
quadratic fields}

\medskip

In this and the next section, we shall study two types of $L$-functions which
are closely related to algebraic number fields. We can still apply our method of
calculation described in \S1 efficiently. However the situation changes drastically.
The repeated application of partial summation does not yield good results beyond
rather limited number of times. Thus our calculation cannot be as accurate as in
the case of modular forms treated in \S2 and \S3.

Let $k$ be a real quadratic field. For simplicity, we assume that the class number of
$k$ is $1$. Let $D$ be the discriminant and $\ep$ be the fundamental unit of $k$.
Let $k_A^{\ti}$ denotes the idele group of $k$. For a finite place $v$ of $k$, let
$k_v$ denote the completion of $k$ at $v$ and $\goth O_v$ denote the ring of
integers of $k_v$. Since $k$ is of class number $1$, we have
$$
k_A^{\ti} = k^{\ti} (\prod _v \goth O_v^{\ti} \ti \RR ^{\ti} \ti \RR ^{\ti} )
\tag{4.1}
$$
where $v$ extends over all finite places of $k$. Let $\chi =\prod _v \chi _v$
be an unramified unitary Hecke character of $k_A^{\ti}$. Let $\s_1$
(resp. $\s _2$) be the identical (resp. non-identical) isomorphism of $k$
into $\RR$ and let $\oo _1$ (resp. $\oo _2$) be the corresponding
archimedean place of $k$. As unitary characters
of $\RR ^{\ti}$, $\chi _{\oo _1}$ and $\chi _{\oo _2}$ take the following form:
$$
\chi _{\oo_j}(x)= \text {sgn} (x)^{m_j} \vert x \vert ^{iv_j} \quad
\text {for} \quad x \in K_{\oo _j}^{\ti} \cong \RR ^{\ti}, \quad j=1,2,
\tag{4.2}
$$
where $m_j=0$ or $1$, $v_j \in \RR$. By (4.1), we see that $\chi$ is completely
determined by $\chi _{\oo _1}$ and $\chi _{\oo _2}$. Since $\chi$ is trivial
on $k^{\ti}$, we must have $\chi (x)=1$ for all
$k^{\ti} \cap (\prod _v \goth O_v^{\ti} \ti \RR ^{\ti} \ti \RR ^{\ti} )$,
which is the group of units of $k$. Therefore we have
$$
(-1)^{m_1+m_2}=1, \qquad
\text {sgn} (\ep ^{\s _1})^{m_1} \text {sgn} (\ep ^{\s _2})^{m_2}
\vert \ep ^{\s _1} \vert ^{iv_1} \vert \ep ^{\s _2} \vert ^{iv_2} =1.
\tag{4.3}
$$
It is easy to see that (4.3) is a necessary and sufficient condition
for $\chi$, which is determined by $\chi _{\oo _1}$ and $\chi _{\oo _2}$,
to be a Hecke character of $k_A^{\ti}$. By (4.3), we have $m_1=m_2$.
Put $m=m_1$. Then (4.3) is equivalent to
$$
\vert \ep \vert ^{i(v_1-v_2)}= \text {sgn} (N(\ep ))^m .
\tag{4.4}
$$
Let $\chi _{\ast}$ be the associated ideal character of $k$. If $(\a )$,
$\a \in k^{\ti}$ is a prime ideal, we have, by definition
$$
\gather
\chi _{\ast}((\a ))=\chi ((\cdots , 1 , \cdots , \a , \cdots , 1 , \cdots )) \\
=\chi (( \a ^{-1}, \cdots , 1 , \cdots , \a ^{-1}, \cdots ))
=\text {sgn} (N(\a ))^m \cdot 
(\vert \a ^{\s _1} \vert ^{iv_1} \vert \a ^{\s _2} \vert ^{iv_2})^{-1} .
\endgather
$$
Here $(\cdots , 1 , \cdots , \a , \cdots , 1 , \cdots ) \in k_A^{\times}$ 
denotes the idele whose $(\alpha)$-component is $\alpha$ and all the other
components are $1$. Hence we have
$$
\chi _{\ast}((\a ))=\text {sgn} (N(\a ))^m \vert \a ^{\s _1} \vert ^{-iv_1}
\vert \a ^{\s _2} \vert ^{-iv_2} \qquad \text {for every} \quad \a \in k^{\ti},
\tag{4.5}
$$
$L(s, \chi )= L(s, \chi _{\ast})= \sum _{(\a )} \chi _{\ast} ((\a ))
N((\a ))^{-s}$. Put
$$
\align
R(s, \chi )&= \vert D \vert ^{s/2} \pi ^{-(s+m)} \pi ^{-i(v_1+v_2)/2}
\Gm ((s+m+iv_1)/2) \\
&\hskip 12.8em \Gm ((s+m+iv_2)/2) L(s, \chi ), \\
R(s, \chi ^{-1})&= \vert D \vert ^{s/2} \pi ^{-(s+m)} \pi ^{i(v_1+v_2)/2}
\Gm ((s+m-iv_1)/2) \\
&\hskip 12.1em \Gm ((s+m-iv_2)/2) L(s, \chi ^{-1}).
\endalign
$$ 
Then the functional equation is (cf. Weil [W], Langlands [LL])
$$
R(s, \chi )= (-1)^m \chi _{\ast}((d))R(1-s, \chi ^{-1}),
\tag{4.6}
$$
where $(d)$ denotes the different of $k$
over $\QQ$. Since $\overline {R(s, \chi )}=R(\bar s, \chi ^{-1})$,
we can put (4.6) in the form of (1.2):
$$
R(1- \bar s, \chi )=(-1)^m \chi _{\ast}((d)) \overline {R(s, \chi )}.
\tag{4.7}
$$
We get
$$
\overline {\chi _{\ast}((d))^{-1/2}R(s, \chi )}=
(-1)^m \chi _{\ast}((d))^{-1/2}R(s, \chi ) \qquad \text {if} \quad
\Re (s)=1/2.
\tag{4.8}
$$
Hence $\chi _{\ast}((d))^{-1/2}R(s, \chi )$ takes real or pure imaginary
values on the critical line according as $m=0$ or $1$. We also note that
we may assume $v_2=0$ without losing any generality since the choice
of $v_2$ can be taken into account as the shift of the variable $s$.
Then, if $m=0$, we have
$v_1=-2n\pi /\log \ep$ with $n \in \ZZ$ by (4.4). We denote this Hecke
character by $\chi _n$. By (4.5), we have
$$
(\chi _n)_{\ast}((\a ))= \vert \a \vert ^{2n\pi i /\log \ep },
\qquad \a \in k^{\ti} .
\tag{4.9}
$$
If $m=1$, by (4.4), we have $v_1=-2n\pi /\log \ep$ or
$v_1=-(2n+1)\pi /\log \ep$ with $n \in \ZZ$  according as
$N(\ep )=1$ or $N(\ep )=-1$.
We denote this Hecke character by $\chi _n^{\pr}$.
By (4.5), we have
$$
(\chi _n^{\pr})_{\ast}((\a ))= 
\cases
\text {sgn} (N(\a )) \vert \a \vert ^{2n\pi i /\log \ep }
\quad \text {if} \quad N(\ep )=1, \\
\text {sgn} (N(\a )) \vert \a \vert ^{(2n+1)\pi i /\log \ep }
\quad \text {if} \quad N(\ep )=-1.
\endcases
\tag{4.10}
$$

As our first example, we take $k=\QQ (\sqrt {2})$. We have 
$\ep =\sqrt {2} +1$, $(d)=(\sqrt {2})^3$, $\vert D \vert =8$. We are
going to study $L(s, \chi _1)$ applying our summation method. For
$s=\frac {1}{2}+it$, $t=15$ and $50$, the values of
$$
R_j=\Re (\chi _{\ast}((d))^{-1/2} \exp (i\vt _0 (t)) S_N^{(j)}), \qquad
I_j=\Im (\chi _{\ast}((d))^{-1/2} \exp (i\vt _0 (t)) S_N^{(j)})
$$
are given in Tables 4.1a and in 4.1b respectively. Here
$\chi _{\ast}((d))^{-1/2}= \sqrt {2}^{-3\pi i/\log \epsilon}$,
$\vt _0(t)= \arg (8^{s/2} \pi ^{-s} \pi ^{2\pi i/\log \epsilon}
\Gamma ((s-2\pi i/\log \epsilon)/2) \Gamma (s/2))$.
From this table, it is evident that $R_j$'s for
higher $j$ do not give good results. We can judge, from the values of
$\vert I_j \vert$, $R_2$ gives the best result, then $R_3$, $R_1$ in this
order. We must be more cautious than in \S2 and \S3 about the accuracy
of the value $e^{i\vt _0 (t)}L(s, \chi _1)$. For example, let $t=15$. 
We empirically judge that $e^{i\vt _0 (t)}L(s, \chi _1)=2.17375$ with error
$\approx 10^{-5}$ from $R_2$ and $I_2$. We have constructed Table 4.2 
in which zeros on the critical line $\Re (s)=1/2$ are listed in the range
$\vert \Im (s) \vert \le 30$. Here, for $n \ge 1$, $t_n$ (resp. $t_{-n}$)
denotes the $n$-th zero of $L(s, \chi _1)$, $s=1/2+it$ on the critical
line for $0 \le t \le 30$ (resp. $0 \le -t \le 30$).

Let us examine the Riemann hypothesis for $L(s, \chi _1)$ in the range
$0 \leq \Im (s) \leq 100$. We have observed $84$ sign changes of
$e^{i\vt _0 (t)} L(\frac {1}{2} +it, \chi _1)$ for $0 \leq t \leq 100$.
For $T>0$, let $N(T)$ denote the number of
zeros of $L(s, \chi _1)$ counted with multiplicity in the domain
$$
0 < \Re (s) < 1, \qquad 0 \le \Im (s) \le T.
$$
Taking $\de =1/2$ in (1.7) , we get
$$
N(T)=\pi ^{-1}\vt (T)+\pi ^{-1} \Delta\text {arg} (L(s, \chi _1)).
\tag{4.11}
$$
Since $L(s, \chi _1) \ne 0$ if $\Re (s) \ge 1$, $\Delta\text {arg} (L(s, \chi _1))$
equals the variation of the argument of $L(s, \chi _1)$ along the line segments
$L_1=[\frac{1}{2}, 1+\mu ]$, $L_2=[1+\mu , 1+ \mu +iT]$,
$L_3=[1+ \mu +iT, \frac{1}{2}+iT]$ for any $\mu >0$. Take $T=100$ and $\mu =1$.
We have $\pi ^{-1}\vt (100)=84.8864 \cdots$. Hence if we can show
$-\pi <  \Delta \text {arg} (L(s, f)) < 0$, then we can conclude that $N(100)=84$.
For $L_1$, we divide it into $15$ intervals. We observed that $L(s, \chi _1)$
moves from $0.3482+0.0712 \sqrt {-1}$ to $0.8011 + 0.0969 \sqrt {-1}$
keeping $\Re (L(s, \chi _1))>0$. For $L_2$, we can show without difficulty that
$\Re (L(s, \chi _1))>0$ on $L_2$. For $L_3$, we divide it into $150$ small intervals. 
We observed that $L(s, \chi _1)$ moves from $0.8159+0.1227 \sqrt {-1}$ to
$-0.0110 -0.0072 \sqrt {-1}$ when $s$ moves from $2+100 \sqrt {-1}$ to
$\frac {1}{2}+100 \sqrt {-1}$; $L(s, \chi _1)$ never crossed the half line
$\Im (L(s, \chi _1))=0$, $\Re (L(s, \chi _1)) \leq 0$. Hence we get $N(100)=84$.
The Riemann hypothesis holds and all zeros of $L(s, \chi _1)$ are simple zeros
in this range.

By (4.11), we should have
$$
N(100)=84.8864 \cdots - \pi ^{-1} \arctan (\frac {712}{3482})
-(1- \pi ^{-1} \arctan (\frac {72}{110}))=84.0067 \cdots .
$$
The error is about $6.7 \times 10^{-3}$ and this is much bigger than the usual
error inherent in our calculations. The reason is that $L(s, \chi _1)$ takes rather small value
at $\frac {1}{2}+100 \sqrt {-1}$ ~; such an error can be made much smaller in the following way.
We take $T=101$. We have observed $86$ sign changes of
$e^{i\vt _0 (t)} L(\frac {1}{2} +it, \chi _1)$ for $0 \leq t \leq 101$.
We have $\pi ^{-1}\vt (101)=86.0881 \cdots$. We divide
$[2+101 \sqrt {-1}, \frac{1}{2}+101 \sqrt {-1}]$ into $150$ small intervals. 
We observed that $L(s, \chi _1)$ moves from $0.9322+0.2563 \sqrt {-1}$ to
$1.7137-0.1288 \sqrt {-1}$. Hence by (4.11), we have
$$
N(101)=86.0881 \cdots - \pi ^{-1} \arctan (\frac {712}{3482})
-\pi ^{-1} \arctan (\frac {1288}{17137})=86
$$
with error less than $10^{-4}$.

In Table 4.3, we have listed $44$ examples of $L(s, \chi _n)$
for which we made experiments in the range $0 \leq \Im (s) \leq T$;
$N(T)$ denotes the number of zeros of $L(s, \chi _n)$ in the domain
$0 < \Re (s) < 1$, $0 \leq \Im (s) \leq T$. We found that all zeros
in the ranges of Table 4.3 lie on the critical line and are simple.

\bigskip

\centerline
{\S5. Artin $L$-functions}

\medskip

Let $k$ be the minimal splitting field of the irreducible polynomial
$f(X)=X^5-X+1$ over $\QQ$. Then $k \supset \QQ (\sqrt {19 \cdot 151})=k_0$,
$k$ is unramified over $k_0$, $\Ga (k/\QQ ) \cong S_5$,
$\Ga (k/k_0) \cong A_5$. The discriminant $\Delta$ of a root of $f(X)$ is
$19 \cdot 151$. This example is due to E. Artin (cf. Lang [LG], p. 121).
Let $\rho$ be an irreducible $4$-dimensional representation of $S_5$
whose character $\chi _{\rho}$ is given as follows.
$$
\spreadmatrixlines{3pt}
\matrix
\text {conjugacy class} & (1) & (12) & (123) & (1234) & (12)(34) & (12)(345) & (12345) \\
\chi _{\rho} & 4 & 2 & 1 & 0 & 0 & -1 & -1
\endmatrix
$$
Since $S_5$ does not have a subgroup of index $4$, $\rho$ is not
monomial. We have $L(s, \rho )= \prod _p L_p(s, \rho )$ for $\Re (s)>1$
with the Euler $p$-factor $L_p(s, \rho )$. We can compute $L_p(s, \rho )$
as follows. First we assume that a prime number $p$ is unramified in $k$,
i.e., $p \ne 19$, $151$. Then we see easily that $L_p(s, \rho )^{-1}$ equals
$$
\alignat 4
(1-p^{-s})^4 &&\qquad &\text {if} \quad \s _p&=&\{(1)\}, \hskip 6.5em\\
(1-p^{-s})^3(1+p^{-s}) &&\qquad &\text {if} \quad \s _p&=&\{(12)\}, \hskip 6.5em\\
(1-p^{-s})^2(1+p^{-s}+p^{-2s}) &&\qquad &\text {if} \quad \s _p&=&\{(123)\}, \hskip 6.5em\\
(1-p^{-s})(1+p^{-s})(1+p^{-2s}) &&\qquad &\text {if} \quad \s _p&=&\{(1234)\}, \hskip 6.5em\\
1+p^{-s}+p^{-2s}+p^{-3s}+p^{-4s} &&\qquad &\text {if} \quad \s _p&=&\{(12345)\},\hskip 6.5em\\
(1-p^{-s})^2(1+p^{-s})^2 &&\qquad &\text {if} \quad \s _p&=&\{(12)(34)\}, \hskip 6.em\\
(1-p^{-s})(1+p^{-s})(1+p^{-s}+p^{-2s}) &&\qquad &\text {if} \quad \s _p&=&\{(12)(345)\} .\hskip 6.5em
\endalignat
$$
Here $\s _p$ denotes the Frobenius conjugacy class of $p$ and $\{ \tau \}$ 
denotes the conjugacy class of $\tau \in S_5$.

Let $p=19$ or $151$. Let $I_{\goth p}$ denote the inertia group of a prime factor
$\goth p$ of $p$ in $k$. By definition, we have
$$
L_p(s, \rho )^{-1} = \det (1- (\rho (\s _{\goth p}) \vert V^{I_{\goth p}}) \cdot p^{-s}) .
\tag{5.1}
$$
Here $V$ denotes the representation space of $\rho$, $V^{I_{\goth p}}$
the subspace of $I_{\goth p}$-fixed vectors and $\s _{\goth p}$ a Frobenius
of $\goth p$ which is determined modulo $I_{\goth p}$. Since $k$ is unramified
over $k_0$, it is obvious that $\vert I_{\goth p} \vert =2$,
$I_{\goth p} \nsubseteq \Ga (k/k_0) \cong A_5$. Hence we may assume that
$I_{\goth p}$ is generated by $(12)$ choosing a suitable $\goth p$ lying over $p$.
Let $D_{\goth p}$ denote the decomposition group of $\goth p$. Then
$D_{\goth p} \triangleright I_{\goth p}$ and $D_{\goth p}/I_{\goth p}$
is generated by $\s _{\goth p} \mod I_{\goth p}$. We have
$$
N_{S_5}(I_{\goth p}) = I_{\goth p} \ti S_3,
$$
where $N_{S_5}(I_{\goth p})$ denotes the normalizer of $I_{\goth p}$ in
$S_5$ and $S_3$ denotes the permutation group on three letters $\{ 3, 4, 5 \}$.

Let $p=19$. Then
$$
f(X) \equiv (X-6)^2(X^3+12X^2+13X+9) \mod 19
$$
is the factorization of $f(X) \mod p$ into irreducible factors in $(\ZZ /p\ZZ )[X]$.
Therefore the residue field extension $\goth O_k/\goth p$ of $\ZZ /p\ZZ$
contains the cubic extension of $\ZZ /p\ZZ$, where $\goth O_k$ denotes
the ring of integers of $k$. Hence we immediately obtain 
$D_{\goth p}/I_{\goth p} \cong \ZZ /3\ZZ$ and that $\s _{\goth p}$ may be taken as 
$(345) \in S_3$ . Now we find easily that
$$
L_p(s, \rho )^{-1} =1-p^{-3s} \qquad \text {if} \quad p=19.
$$

Let $p=151$. Then
$$
f(X) \equiv (X-39)^2(X-9)(X^2+87X+61) \mod 151
$$
is the factorization of $f(X) \mod p$ into irreducible factors in $(\ZZ /p\ZZ )[X]$.
By a similar consideration as above, we find that
$D_{\goth p}/I_{\goth p} \cong \ZZ /2\ZZ$ and that $\s _{\goth p}$ may be taken
as  $(34) \in S_3$. We obtain
$$
L_p(s, \rho )^{-1} =(1+p^{-s})(1-p^{-s})^2 \qquad \text {if} \quad p=151.
$$

Let $f(\rho )$ denote the Artin conductor of $\rho$. We easily obtain
$$
f(\rho )=19 \cdot 151.
$$
For example, let $p=19$ and $\goth p$ be as above. We have shown
$\Ga (k_{\goth p}/\QQ _{19}) \cong D_{\goth p} \cong \ZZ /2\ZZ \ti \ZZ /3\ZZ$.
We find that the restriction of $\rho$ to $D_{\goth p}$ splits into a direct sum
of four one dimensional representations of $D_{\goth p}$ such that three of them are
unramified and one is ramified. Hence the exponent of $19$ in $f(\rho )$ is $1$.

We take an isomorphism $\s$ of $k$ into $\CC$ and let $c \in \Ga (k/\QQ ) \cong S_5$
be the restriction of the complex conjugation to $k$. Then 
$c \in \Ga (k/k_0) \cong A_5$. Hence $c$ is conjugate to $(12)(34)$ in $S_5$.
Let $\Ga (\CC /\RR )$ be identified with the decomposition group
$\la c \ra$ of the archimedean place of $k$ which corresponds to $\s$.
The restriction of $\rho$ to $\Ga (\CC /\RR )$ splits into a direct sum of
two trivial representations and two non-trivial representations. Therefore
the Gamma factor to go with $L(s, \rho )$ is given by (cf. Langlands [LL])
$$
(\pi ^{-s/2} \Gm (s/2))^2 (\pi ^{-(s+1)/2} \Gm ((s+1)/2))^2 .
$$
Put
$$
R(s, \rho )= (19 \cdot 151)^{s/2} \pi ^{-2s} \Gm (s/2)^2
\Gm ((s+1)/2)^2 L(s, \rho ).
$$
Since $\rho$ is equivalent to its contragredient, we have the functional equation
$$
R(s, \rho )= \kappa R(1-s, \rho ),
\tag{5.2}
$$
where $\kappa =\pm 1$ is the Artin root number attached to $\rho$. Let
$\psi$ be the additive character of $\QQ _A/\QQ$ such that
$$
\align
\psi _{\oo} (x)&= \exp (2\pi \sqrt {-1}x), \qquad x \in \QQ _{\oo} \cong \RR , \\
\psi _p(x)&=\exp (-2\pi \sqrt {-1} \,  \text {Fr}(x)), \qquad x \in \QQ _p ,
\endalign
$$
where $\text {Fr}$ denotes the fractional part of $x$. By a theorem of Langlands,
we have
$$
\kappa = \prod _v \ep (\frac {1}{2}, \rho _v, \psi _v)
$$
with the $\ep$-factor defined in [LL]. By the above considerations, 
we easily get
$$
\align
\ep (\frac {1}{2}, \rho _p, \psi _p)&=
\cases
1 \qquad \text {if} \quad p \ne 19, 151, \\
i \qquad \text {if} \quad p=19 \ \text {or} \ 151, 
\endcases \\
\ep (\frac {1}{2}, \rho _{\oo}, \psi _{\oo})&=i^2. 
\endalign
$$
Hence we obtain
$$
\kappa =1.
\tag{5.3}
$$

The values of
$$
R_j=\Re (\exp (i\vt (t)) S_N^{(j)}), \qquad
I_j=\Im (\exp (i\vt (t)) S_N^{(j)})
$$
for $L(s, \rho )$, $s=\frac {1}{2}+it$, $t=5$ are given in Table 5.1. Here
$$
\vt (t)= \text {arg} ((19 \cdot 151)^{s/2} \pi ^{-2s} 
\Gm (s/2)^2\Gm ((s+1)/2)^2 ), \qquad s=\frac {1}{2}+it.
$$
In Table 5.2, we give the values of $u_n$, the $n$-th zero of
$L(s, \rho )$, $s=\frac{1}{2}+iu$ on the critical line for $0 \leq u \leq 10$.

\bigskip

\centerline {\S 6. Estimation of errors in our calculations}

\medskip

The most serious defect of our method of calculation is that we do not have
rigorous controle of error estimates. In previous sections, we regarded the magnitude
of $\Im (e^{i\vt (t)} S_N^{(l)})$  (resp. $\Re (e^{i\vt (t)} S_N^{(l)})$ )
 as a rough measure of errors from the true value, when $e^{i\vt (t)} L(\sigma +it)$ 
should be real (resp. pure imaginary). In this section, we shall present several data
which support this practice.

Suppose that the functional equation (1.2) for $L(s)$ holds. Then we have 
\newline $e^{i\vt (t)} L(k/2 +it) \in \bold R$, $t \in \bold R$ where 
$\vt (t)=\arg (\kappa _1 N^s \prod _{i=1}^m \Gm (b_is+c_i))$,
$s=k/2+it$ in the notation of (1.4).
Take $0<t_1<t_2$ so that $e^{i\vt (t)}$ rotates on
the unit circle exactly once when $t$ moves from $t_1$ to $t_2$. We expect that
$\max _{t_1 \leq t \leq t_2} \vert \Im (e^{i\vt (t)} S_N^{(l)}) \vert$
can be used as the measure of errors. More explicitly, it seems plausible that
$$
\vert \Re (e^{i\vt (t)} S_N^{(l)}) - L(k/2 +it) \vert \leq
10 \max _{t_1 \leq t \leq t_2} \vert \Im (e^{i\vt (t)} S_N^{(l)}) \vert .
\tag{6.1}
$$
In examples below, we use $\Re (S_M^{(p)})$ as a substitute for $L(k/2 +it)$ 
taking large $M$ and $p$ (except for in Example 6),
and examine the ratio of two terms in (6.1) for $S_N^{(l)}$ taking relatively small
$N$ and $l$. The results are given in Table 6.1.

\medskip

Example 1.  We take the primitive form $f \in S_8(\Gamma _0(2))$ and consider
the $L$-function $L(s,f)$, $s=4+it$. When $t$ moves from $97.9$ to $100$,
$e^{i\vt (t)}$ rotates on the unit circle approximately once. We calculated the ratio
$$
r_{5N}^{(1)}=\frac 
{\max _{97.9 \leq t \leq 100} ( \vert \Re (e^{i\vt (t)}S_{2000}^{(5N)} -
e^{i\vt (t)}S_{10000}^{(35)}) \vert )}
{\max _{97.9 \leq t \leq 100}  (\vert \Im (e^{i\vt (t)}S_{2000}^{(5N)}) \vert )}
$$
for $1 \leq N \leq 6$, dividing $[97.9, 100]$ into $21$ intervals of length $0.1$.

\medskip

Example 2.  We take $g \in S_{9/2}(\Gamma _0(4))$ as in \S 3. When $t$ moves from $98.2$ to $100$,
$e^{i\vt (t)}$ rotates on the unit circle approximately once. We calculated the ratio
$$
r_{5N}^{(2)}=\frac 
{\max _{98.2 \leq t \leq 100} ( \vert \Re (e^{i\vt (t)}S_{2000}^{(5N)} -
e^{i\vt (t)}S_{10000}^{(35)}) \vert )}
{\max _{98.2 \leq t \leq 100}  (\vert \Im (e^{i\vt (t)}S_{2000}^{(5N)}) \vert )}
$$
for $1 \leq N \leq 6$, dividing $[98.2, 100]$ into $18$ intervals of length $0.1$.

\medskip

Example 3. We take $\Delta \in S_{12}(SL(2, \bold Z))$ and consider the $L$-function
$L^{(3)}(s, \Delta )$, $s=17+it$, attached to the third symmetric power representation
of $GL(2)$. When $t$ moves from $17.5$ to $20$,
$e^{i\vt (t)}$ rotates on the unit circle approximately once. We calculated the ratio
$$
r_{5N}^{(3)}=\frac 
{\max _{17.5 \leq t \leq 20} ( \vert \Im (e^{i\vt (t)}S_{2000}^{(5N)} -
e^{i\vt (t)}S_{10000}^{(35)}) \vert )}
{\max _{17.5 \leq t \leq 20}  (\vert \Re (e^{i\vt (t)}S_{2000}^{(5N)}) \vert )}
$$
for $1 \leq N \leq 6$, dividing $[17.5, 20]$ into $25$ intervals of length $0.1$.
In this example, we have normalized $\vt (t)$ as in \S2 so that
$e^{i\vt (t)}L^{(3)}(17+it, \Delta )$ is pure imaginary.

\medskip

Example 4. We consider the $L$-function
$L^{(4)}(s, \Delta )$, $s=45/2+it$, attached to the fourth symmetric power representation
of $GL(2)$. When $t$ moves from $7.2$ to $10$,
$e^{i\vt (t)}$ rotates on the unit circle approximately once. We calculated the ratio
$$
r_{5N}^{(4)}=\frac 
{\max _{7.2 \leq t \leq 10} ( \vert \Re (e^{i\vt (t)}S_{2000}^{(5N)} -
e^{i\vt (t)}S_{10000}^{(35)}) \vert )}
{\max _{7.2 \leq t \leq 10}  (\vert \Im (e^{i\vt (t)}S_{2000}^{(5N)}) \vert )}
$$
for $1 \leq N \leq 6$, dividing $[7.2, 10]$ into $28$ intervals of length $0.1$.

\medskip

Example 5. We consider the Artin $L$-function treated in \S 5.
When $t$ moves from $4$ to $5.8$,
$e^{i\vt (t)}$ rotates on the unit circle approximately once. We calculated the ratio
$$
r_{2N}^{(5)}=\frac 
{\max _{4 \leq t \leq 5.8} ( \vert \Re (e^{i\vt (t)}S_{10000}^{(2N)} -
e^{i\vt (t)}S_{100000}^{(10)}) \vert )}
{\max _{4 \leq t \leq 5.8}  (\vert \Im (e^{i\vt (t)}S_{10000}^{(2N)}) \vert )}
$$
for $0 \leq N \leq 5$, dividing $[4, 5.8]$ into $18$ intervals of length $0.1$.

\medskip

Example 6. We consider the Hecke $L$-function $L(s, \chi _1)$ 
for $k= \bold Q( \sqrt 2)$ treated in \S 4.
When $t$ moves from $47.9$ to $50$,
$e^{i\vt (t)}$ rotates on the unit circle approximately once. We calculated the ratio
$$
r_N^{(6)}=\frac 
{\max _{47.9 \leq t \leq 50} ( \vert \Re (e^{i\vt (t)}S_{10000}^{(N)} -
e^{i\vt (t)}S_{100000}^{(2)}) \vert )}
{\max _{47.9 \leq t \leq 50}  (\vert \Im (e^{i\vt (t)}S_{10000}^{(N)}) \vert )}
$$
for $0 \leq N \leq 5$, dividing $[47.9, 50]$ into $21$ intervals of length $0.1$.

\bigskip

\centerline {\S7. A comparison with the explicit formula}

\medskip

We take the new form $f \in S_8(\Gamma _0(2))$ treated in \S3. Let
$\pi$ be the irreducible unitary automorphic representation of
$GL(2, \bold Q_A)$ which corresponds to $f$. We have
$$
L_f(s, \pi )=L(s + \frac {7}{2}, f),
$$
where $L_f(s, \pi )$ denotes the finite part of the Jacquet-Langlands
$L$-function attached to $\pi$. Let
$$
L_f(s, \pi )= \prod _p \left\lbrack 
(1 - \alpha _p p^{-s})(1 - \beta _p p^{-s}) \right\rbrack ^{-1}
$$
be the Euler product of $L_f(s, \pi )$. For $p=2$, the Euler
$2$-factor degenerates so that $\alpha _2= -1/\sqrt {2}$,
$\beta _2=0$. For $p \neq 2$, we have
$\vert \alpha _p \vert = \vert \beta _p \vert =1$
by the Ramanujan-Petersson conjecture proved by P. Deligne. We have
the explicit formula
$$
\aligned
\sum _p \sumprime _{1 \leq n, p^n \leq x}
(\alpha _p^n + \beta _p^n) \log p =
&- \lim _{T \rightarrow +\infty} \sum _{\vert \Im (\rho ) \vert <T}
\frac {x^{\rho}} {\rho} - 
\frac {L^{\prime} _f(0, \pi )} {L_f(0, \pi ) } \\
& + \log (\frac {\sqrt {x} +1}{\sqrt {x} -1})
-2x^{-1/2} -\frac {2}{3}x^{-3/2} -\frac {2}{5} x^{-5/2}
\endaligned
\tag{7.1}
$$
for $x>1$. Here $\sumprime$ means that the term
$(\alpha _p^n + \beta _p^n) \log p$ should be multiplied
by $1/2$ when $p^n=x$; $\rho$ extends over zeros of $L_f(s, \pi )$
such that $0< \Re (\rho )<1$. This formula can be shown in
the usual way as in Ingham [In], p. 77--80. The last term
$$
g(x):= \log (\frac {\sqrt {x} +1}{\sqrt {x} -1})
-2x^{-1/2} -\frac {2}{3}x^{-3/2} -\frac {2}{5} x^{-5/2},
$$
which is equal to
$$
2 \sum _{k=3} ^{\infty} \frac {x^{-(2k+1)/2}} {2k+1},
$$
represents the contribution of the trivial zeros of $L_f(s, \pi )$;
they are at $s=-\frac {7}{2}$, $-\frac {9}{2}$, $-\frac {11}{2}$,
$\cdots$.

Now it seems very interesting to compare both sides of (7.1) numerically
using the zeros of $L(s,f)$ given in Table 3.3. We approximate
$$
\lim _{T \rightarrow +\infty} \sum _{\vert \Im (\rho ) \vert <T}
\frac {x^{\rho}} {\rho}
$$
by
$$
h_{69}(x):= \sum _{n=1}^{69} \sqrt {x}
(t_n^2 +\frac {1}{4} )^{-1}
\lbrack \cos (t_n \log x)+ 2t_n \sin (t_n \log x) \rbrack
\tag{7.2}
$$
with $t_n$ given in Table 3.3. We have
$$
\align
L_f(0, \pi )&=L(\frac {7}{2}, f)=0.5942254156 \cdots , \quad
L_f^{\prime} (0, \pi )=L^{\prime} (\frac {7}{2}, f)=
0.1875716234 \cdots , \\
\frac {L_f^{\prime} (0, \pi )}{L_f(0, \pi )}&=0.3156573558 \cdots .
\endalign
$$
To obtain these values, we simply applied our repeated abel
summation technique as before, though more rigorous evaluation
could be made in this case.

In Figure 7.1, we drawed the graphs of the \lq\lq step function \hskip -0.3em \rq\rq
$$
\sum _p \sumprime _{1 \leq n, p^n \leq x}
(\alpha _p^n + \beta _p^n) \log p
$$
and $-h_{69}(x)+g(x)-0.3156$ for $1.1 \leq x \leq 20$. (We have
used \lq\lq Mathematica \hskip -0.3em \rq\rq \, \, to make Figure 7.1.)
The coincidence seems fine.

\bigskip

\centerline {\S 8. Sample programs}

\medskip

In this section, we shall present a few sample programs to compute $L(s,f)$,
$f \in S_8(\Gamma _0(2))$ (cf. \S 3). All programs, which are ready to be
executed, are written in UBASIC created by Y. Kida. Let
$$
f(z)=(\eta (z)\eta (2z))^8=\sum _{n=1}^{\oo} a_nq^n .
$$
Using Program A, we can compute $a_n$ for $1 \leq n \leq M$ for any
$M$, $1 \leq M \leq 10^4$. From line $50$ to $150$, the coefficients $A(n)$
in $\eta (z)= q^{1/24} \sum _{n=1}^{\infty} A(n)q^{n-1}$ are computed for
$1 \leq n \leq M$ using Euler's formula
$$
\eta (z)= q^{1/24}\prod _{n=1}^{\infty} (1-q^n)
=q^{1/24} \sum _{n=-\infty}^{\infty} (-1)^n q^{n(3n+1)/2}.
$$
From line $160$ to $220$, the coefficients $B(n)$ in
$\eta (z)\eta (2z)=q^{1/8}\sum _{n=1}^{\infty} B(n)q^{n-1}$ are computed
for $1 \leq n \leq M$ by $B(I)= \sum _{J+2L=I+2} A(J)A(L)$.
From line $270$ to $390$, the expansion of $\eta (z)\eta (2z)$ is raised to
the eighth power; the final result will be stored in the data file \lq\lq wt8 \rq\rq .

Program B computes the value of $e^{i\vt _f (t)}L(4+it,f)$ for $t=100$.
To save the space, this program gives the values of $e^{i\vt _f (t)}S_N^{(l)}$
for $N=2000$, $0 \leq l \leq 40$. By point $15$ command (line $20$),
UBASIC gives the precision to the $70$-th digit.
From line $300$ to $580$, the value of
$e^{i\vt _f (t)}$ will be computed and stored in the variable THE. 
The calculation proceeds as follows. We have
$$
\vt _f (t)= \arg (2^{s/2}(2\pi )^{-s} \Gamma (s))_{s=4+it}
= (\frac {1}{2} \log 2 - \log 2\pi )t +\arg \Gamma (4+it).
$$
We also have
$$
\arg \Gamma (z)= \arg \Gamma (z+1)- \arctan (\Im (z)/ \Re (z)),
\quad \Re (z)>0, \quad \arg \Gamma (z)= \Im (\log \Gamma (z)).
$$
By these formulas, it suffices to compute $\log \Gamma (z+100)$ for
$z=4+it$. We have the asymptotic expansion (cf. [WW], p. 252)
$$
\log \Gamma (z) \sim (z- \frac {1}{2}) \log z -z + \frac {1}{2} \log 2\pi
+ \sum _{r=1}^{\infty} \frac {(-1)^{r-1} B_r}{2r(2r-1)z^{2r-1}} ,
\tag{8.1}
$$
where $B_r$ denotes the $r$-th Bernoulli number.
In (8.1), we use the terms up to $r=10$. We can get an approximation
of $e^{i\vt _f (t)}$ which is accurate at least to the $40$-th digit (cf. [WW], p. 252).

From line $660$ to $740$, the values of
$$
S_N^{(0)}=\sum _{n=1}^N a_nn^{-s}, \quad
s_n^{(l)}=\sum _{n=1}^{N} s_n^{(l-1)}, \quad 1 \leq l \leq 40,
\quad (s_n^{(0)}=a_n)
$$
are computed. From line $770$ to $1020$, the values of
$S_N^{(l)}$ are computed for $1 \leq l \leq 40$ using (1.10) and (1.11).
We use the approximation
$$
u_N^{(l)}=N^{-s} \sum _{k=l}^{100}
(\sum _{m=1}^l (-1)^m \binom lm m^k)
(-1)^k \frac {s(s+1) \cdots (s+k-1)} {k !} N^{-k}.
$$
The error from the truncation by $100$ is negligible (cf. \S 1). The variable $U(K)$ stands for
$(-1)^k \frac {s(s+1) \cdots (s+k-1)} {k !}$;
the variable $Co$ stands for 
$\sum _{m=1}^l (-1)^m \binom lm m^k$;
the variable $T(I)$ stands for 
$$
\sum _{k=l}^{100}
(\sum _{m=1}^l (-1)^m \binom lm m^k)
(-1)^k \frac {s(s+1) \cdots (s+k-1)} {k !} N^{-k}.
$$
In the line $1000$, $X=S_N^{(I)}$ is multiplied by $e^{i\vt _f (t)}$.

\newpage

\centerline {Program A}

\medskip

\item {10} word 8
\item {20} point 2
\item {30} dim A(10000),B(10000)
\item {40} input M
\item {50} A1=sqrt(24*M+1)
\item {60} M1=int((A1+1)/6)
\item {70} A(1)=1
\item {80} for I=1 to M1
\item {90} I1=I$-$2$\ast$int(I/2)
\item {100} I2=1$-$2$\ast$I1
\item {110} J=int((3$\ast$I$\ast$I+I)/2)+1
\item {120} A(J)=A(J)+I2
\item {130} J=int((3$\ast$I$\ast$I$-$I)/2)+1
\item {140} A(J)=A(J)+I2
\item {150} next I
\item {160} for I=1 to M
\item {170} I1=int((I+1)/2)
\item {180} for L=1 to I1
\item {190} J=I$-$2$\ast$L+2
\item {200} B(I)=B(I)+A(J)$\ast$A(L)
\item {210} next L
\item {220} next I
\item {230} for I=1 to M
\item {240} A(I)=B(I)
\item {250} B(I)=0
\item {260} next I
\item {270} K=1
\item {280} for I=1 to M
\item {290} print K,I
\item {300} for J=1 to I
\item {310} B(I)=B(I)+A(J)$\ast$A(I+1$-$J)
\item {320} next J
\item {330} next I
\item {340} for I=1 to M
\item {350} A(I)=B(I)
\item {360} B(I)=0
\item {370} next I
\item {380} K=K+1
\item {390} if K$<$4 then goto 280
\item {400} open "wt8" for output as $\sharp$1
\item {410} for I=1 to M
\item {420} print I,A(I)
\item {430} print $\sharp$1,A(I)
\item {440} next I
\item {450} close $\sharp$1
\item {460}  end

\bigskip

\centerline {Program B}   

\medskip

\item {10} word 70
\item {20} point 15
\item {30} dim Bn(20),Bd(20),C(2000),T(100),U(200),Sm(50)
\item {40} Ab=40
\item {50} Bn(1)=1
\item {60} Bd(1)=6
\item {70} Bn(2)=1
\item {80} Bd(2)=30
\item {90} Bn(3)=1
\item {100} Bd(3)=42
\item {110} Bn(4)=1
\item {120} Bd(4)=30
\item {130} Bn(5)=5
\item {140} Bd(5)=66
\item {150} Bn(6)=691
\item {160} Bd(6)=2730
\item {170} Bn(7)=7
\item {180} Bd(7)=6
\item {190} Bn(8)=3617
\item {200} Bd(8)=510
\item {210} Bn(9)=43867
\item {220} Bd(9)=798
\item {230} Bn(10)=174611
\item {240} Bd(10)=330
\item {250} A=1/sqrt(3)
\item {260} P=6$\ast$atan(A)
\item {270} S=4+100$\ast$$\sharp$i
\item {280} Ss=S
\item {290} T=im(S)
\item {300} Th=T$\ast$(log(2)/2$-$log(2$\ast$P))
\item {310} U=0
\item {320} for I=1 to 100
\item {330} S1=re(S)
\item {340} S2=im(S)
\item {350} if S1$>$S2 goto 380
\item {360} U=U$-$(P/2)+atan(S1/S2)
\item {370} goto 390
\item {380} U=U$-$atan(S2/S1)
\item {390} S=S+1
\item {400} next I
\item {410} Th=Th+U
\item {420} T=im(S)
\item {430} R=re(S)
\item {440} X1=atan(T/R)
\item {450} X2=(log(R$\ast$R+T$\ast$T))/2
\item {460} X3=X2+X1$\ast$$\sharp$i
\item {470} X4=(S$-$(1/2))$\ast$X3
\item {480} Th=Th+im(X4)$-$T
\item {490} K=10
\item {500} S1=S
\item {510} for I=1 to K
\item {520} X1=Bn(I)/(2$\ast$I$\ast$(2$\ast$I$-$1)$\ast$Bd(I)$\ast$S1)
\item {530} Th=Th+im(X1)
\item {540} S1=$-$S1$\ast$S$\ast$S
\item {550} next I
\item {560} X1=int(Th/(2$\ast$P))
\item {570} X2=Th$-$2$\ast$P$\ast$X1
\item {580} The=exp(X2$\ast$$\sharp$i)
\item {590} S=Ss
\item {600} X=0
\item {610} open "wt8" for input as $\sharp$1
\item {620} for M=1 to 2000
\item {630} input $\sharp$1,C(M)
\item {640} next M
\item {650} close $\sharp$1
\item {660} for M=1 to 2000
\item {670} Sm(1)=Sm(1)+C(M)
\item {680} for I=2 to Ab
\item {690} Sm(I)=Sm(I)+Sm(I$-$1)
\item {700} next I
\item {710} X1=log(M)
\item {720} X2=exp($-$X1$\ast$S)
\item {730} X=X+C(M)$\ast$X2
\item {740} next M
\item {750} Z=The$\ast$X
\item {760} print Z
\item {770} N=2001
\item {780} A=log(N)
\item {790} N1=exp($-$A$\ast$S)
\item {800} U(1)=$-$S
\item {810} for K=2 to 100
\item {820} U(K)=$-$U(K$-$1)$\ast$(S+K$-$1)/K
\item {830} next K
\item {840} for I=1 to Ab
\item {850} T(I)=0
\item {860} if I$>$1 then goto 890
\item {870} X=X$-$Sm(1)$\ast$N1
\item {880} goto 1000
\item {890} Ii=I$-$1
\item {900} for K=Ii to 100
\item {910} Co=$-$Ii
\item {920} A1=$-$Ii
\item {930} for L=2 to Ii
\item {940} A1=$-$A1$\ast$(Ii$-$L+1)/L
\item {950} Co=Co+(A1$\ast$(L$\wedge$K))
\item {960} next L
\item {970} T(I)=T(I)+Co$\ast$U(K)/(N$\wedge$K)
\item {980} next K
\item {990} X=X$-$Sm(I)$\ast$N1$\ast$T(I)
\item {1000} Z=The$\ast$X
\item {1010} print I,Ss,Z
\item {1020} next I
\item {1030} end

\bigskip

\centerline {\S9. Conjectures}

\medskip

In this section, we shall discuss a few conjectures which emerged
during the process of experiments: No non-trivial coincidences of
zeros of two $L$-functions attached to two non-equivalent irreducible
$\lambda$-adic representations of
$\text {Gal} (\overline {\bold Q}/\bold Q)$
are found. We shall use the framework of automorphic representations
of $GL(n, \bold Q_A)$ to formulate this fact in more general case.

Let $\pi$ be an irreducible unitary cuspidal automorphic representation of 
\newline $GL(n, \bold Q_A)$. The contragredient representation $\tilde \pi$
to $\pi$ is equivalent to the complex conjugate representation
$\bar \pi$ of $\pi$ and we have the functional equations:
$$
L(s, \pi )=\epsilon (s, \pi )L(1-s, \bar \pi ), \qquad
L(s, \bar \pi )=\overline {L(\bar s, \pi )}.
$$
Let $\omega _{\pi}$ be the central character of $\pi$. For $s \in \bold C$,
set $\nu ^s(x)= \vert x \vert _A^s$, $x \in \bold Q _A^{\times}$ where
$\vert x \vert _A$ denotes the idele norm of $x$. We can find a 
$t \in \bold R$ so that $\omega _{\pi}\nu ^{it}$ is a character of
$\bold Q_A^{\times}$ of finite order. Since $\omega _{\pi}\nu ^{it}$
is the central character of $\pi \otimes (\nu ^{it/n} \circ \det )$ and
$L(s+ \frac {it}{n}, \pi )=L(s, \pi \otimes (\nu ^{it/n} \circ \det ))$,
we may assume, without losing substantial generality, that $\omega _{\pi}$
is of finite order.

\proclaim {Conjecture 9.1}  Let $\pi _1$ and $\pi _2$ be irreducible
unitary cuspidal automorphic representations of $GL(n, \bold Q_A)$
with the central characters $\omega _{\pi _1}$ and  $\omega _{\pi _2}$.
We assume that $\pi _1$ is not equivalent to $\pi _2$ and that
$\omega _{\pi _1}$ and $\omega _{\pi _2}$ are of finite order.
Then $L(s, \pi _1)$ and $L(s, \pi _2)$ have no common zeros
in the critical strip $0 < \Re (s) <1$ except for $s=1/2$.
\endproclaim

Remark. If we replace $\bold Q$ by an algebraic number field, 
the assertion is obviously false.

\medskip

As a variant of $9.1$, we can formulate a conjecture on
$L$-functions of motives. Let $E$ be an algebraic number field of
finite degree. Let $M_1$ and $M_2$ be motives over $\bold Q$ 
with coefficients in $E$ of pure weights $w_1$ and $w_2$, of ranks
$n_1$ and $n_2$ respectively. We assume $w_1=w_2$ and put $w=w_1$. 
Fix an embedding $\sigma$ of $E$ into $\bold C$ and let $L(s, M_1)$
(resp. $L(s, M_2)$) be the $L$-function of $M_1$ (resp. $M_2$) with
respect to $\sigma$. For a finite place $\lambda$ of $E$, let
$\rho _i : \text {Gal} (\overline {\bold Q}/\bold Q)
\longrightarrow GL(n_i, E_{\lambda})$, $i=1$, $2$ be the
$\lambda$-adic representation obtained from the $\lambda$-adic
realization of $M_i$. We assume meromorphic continuation of
$L(s, M_i)$, $i=1$, $2$ to the whole complex plane.

\proclaim {Conjecture 9.2} Assume that $\rho _i$, $i=1$, $2$ are
absolutely irreducible and that $\rho _1$ is not equivalent to 
$\rho _2$ for a finite place $\lambda$ of $E$. Then $L(s, M_1)$
and $L(s, M_2)$ have no common zeros in the critical strip
$\frac {w}{2} < \Re (s) < \frac {w}{2}+1$ except for $s=(w+1)/2$.
\endproclaim

Remark. We understand that a pole of order $k$ is a zero of
order $-k$. Conjecture 9.2 implies the (usual) Artin conjecture
except for $s=1/2$.

Remark. The condition on $\rho _i$ implies that $M_1$ and $M_2$
are simple motives with coefficients in $\overline {\bold Q}$.
D. Blasius observed that the analogous conjecture for motives over
a finite field is true under the Tate conjecture (cf. Milne [M], p. 415,
Proposition 2.6).

\newpage

$$
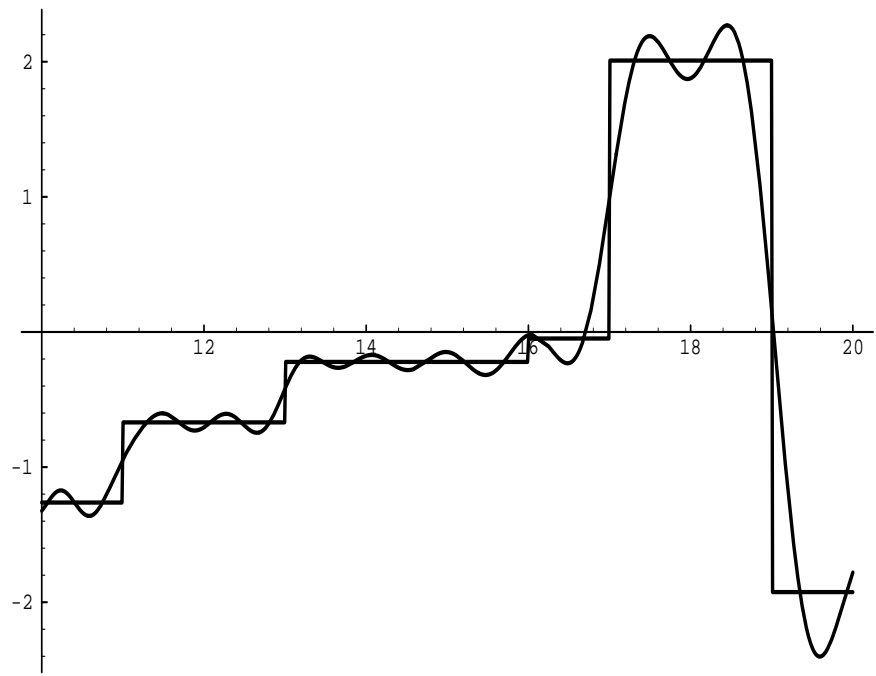
\text{Table 2.1}
$$

\define\qc{\hskip 20pt}
\halign{$\hfil#\hfil$\quad &$\hfil#\hfil$\quad &$\hfil#$\quad
&$\hfil#\hfil$\quad &$\hfil#$\quad &$\hfil#\hfil$\quad
&$\hfil#$ \cr
\noalign{\hrule}
\noalign{\smallskip}
N & I_0 & R_0 \qc & I_5 & R_5 \qc & I_{10} & R_{10} \qc \cr
\noalign{\smallskip}
\noalign{\hrule}
\noalign{\smallskip}
1000& 3.95& -9.3\times 10^{-2}&
3.9415443& 1.3\times 10^{-3}&
3.9237277& 3.0\times 10^{-3} \cr
2000& 4.02& 3.3\times 10^{-2}&
3.9248413& 4.6\times 10^{-3}&
3.9238190& -1.1\times 10^{-5} \cr
4000& 3.90& 1.0\times 10^{-1}&
3.9233635& -2.1\times 10^{-4}&
3.9242209& -1.8\times 10^{-5} \cr
6000& 3.88& -3.9\times 10^{-3}&
3.9239218& -3.6\times 10^{-4}&
3.9242178& -3.7\times 10^{-6} \cr
8000& 3.85& -1.4\times 10^{-2}&
3.9245956& 6.6\times 10^{-5}&
3.9241989& 5.7\times 10^{-6} \cr
10000& 3.92& 2.5\times 10^{-2}&
3.9240769& -8.5\times 10^{-7}&
3.9242054& -3.2\times 10^{-7} \cr
\noalign{\smallskip}
\noalign{\hrule}
}

\medskip

\redefine\qc{\hskip 23pt}
\define\qd{\hskip4.5pt}
\halign{$\hfil#\hfil$\quad &$\hfil#\hfil$\quad &$\hfil#$\quad
&$\hfil#\hfil$\quad &$\hfil#$ \cr
\noalign{\hrule}
\noalign{\smallskip}
N & I_{15} & R_{15} \qc & I_{20} & R_{20} \qc \cr
\noalign{\smallskip}
\noalign{\hrule}
\noalign{\smallskip}
1000& 3.92361007& -1.3\times 10^{-4}&
3.924178604& -4.1\times 10^{-4} \cr
2000& 3.92423322& -3.8\times 10^{-5}&
3.924207205& 4.9 \times 10^{-6} \cr
4000& 3.92420427& 1.8\times 10^{-6}&
3.924203509& -1.8\times 10^{-7} \cr
6000& 3.92420338& 4.6\times 10^{-7}&
3.924203748& -2.9\times 10^{-8} \cr
8000& 3.92420375& -1.9\times 10^{-7}&
3.924203739& 7.2\times 10^{-9} \cr
10000& 3.92420370& 2.4\times 10^{-8}&
3.924203738& -1.0\times 10^{-9} \cr
\noalign{\smallskip}
\noalign{\hrule}
}

\medskip

\redefine\qc{\hskip 23pt}
\redefine\qd{\hskip4.5pt}
\halign{$\hfil#\hfil$\quad &$\hfil#\hfil$\quad &$\hfil#$\quad
&$\hfil#\hfil$\quad &$\hfil#$ \cr
\noalign{\hrule}
\noalign{\smallskip}
N & I_{25} & R_{25} \qc & I_{30} & R_{30} \qc \cr
\noalign{\smallskip}
\noalign{\hrule}
\noalign{\smallskip}
1000& 3.92438069843& 8.8\times 10^{-6} \qd &
3.92417001558& 2.0\times 10^{-4} \qd \cr
2000& 3.92419997465& 1.2\times 10^{-6} \qd &
3.92420393248& -9.9\times 10^{-7} \qd \cr
4000& 3.92420378031& 2.6\times 10^{-8} \qd &
3.92420373301& 3.0\times 10^{-9} \qd \cr
6000& 3.92420373864& 1.5\times 10^{-9} \qd &
3.92420373759& -1.7\times 10^{-10} \cr
8000& 3.92420373800& -2.7\times 10^{-10}&
3.92420373818& -8.9\times 10^{-13} \cr
10000& 3.92420373814& 5.8\times 10^{-11}&
3.92420373812& -1.0\times 10^{-11} \cr
\noalign{\smallskip}
\noalign{\hrule}
}

\medskip

\redefine\qc{\hskip 23pt}
\redefine\qd{\hskip4.5pt}
\halign{$\hfil#\hfil$\quad &$\hfil#\hfil$\quad &$\hfil#$ \cr
\noalign{\hrule}
\noalign{\smallskip}
N & I_{35} & R_{35} \qc \cr
\noalign{\smallskip}
\noalign{\hrule}
\noalign{\smallskip}
1000& 3.923939180351& 4.8\times 10^{-5} \qd \cr
2000& 3.924204540724& 8.0\times 10^{-7} \qd \cr
4000& 3.924203740852& -3.7\times 10^{-9} \qd \cr
6000& 3.924203738244& 7.9\times 10^{-11} \cr
8000& 3.924203738121& -4.4\times 10^{-12} \cr
10000& 3.924203738135& 2.1\times 10^{-12} \cr
\noalign{\smallskip}
\noalign{\hrule}
}

\newpage

$$
\text{Table 2.2}
$$

$$\vbox{\offinterlineskip
\define\qa{\hskip6pt}
\define\qaa{\hskip7pt}
\define\qb{\hskip2.3pt}
\redefine\qc{\hskip28pt}
\redefine\qd{\hskip5.3pt}
\define\vsp{height2pt & \omit & & \omit & & \omit & & \omit & & \omit
& & \omit & \cr}
\define\vspa{height3pt & \omit & & \omit & & \omit & & \omit & & \omit
& & \omit & \cr}
\halign{&\vrule# &\strut\qa\hfil#\hfil\qb &\vrule#
&\strut\qaa#\hfill\qb \cr
\noalign{\hrule}
\vsp
& n & & \qc $t_n$ & & n & & \qc $t_n$ & & n & & \qc $t_n$ & \cr
\vsp
\noalign{\hrule}
\vspa
& 1 & & \qd 0 & & 2 & & \qd 4.1558656464 & & 3 & & \qd 5.5491219562 & \cr
& 4 & & \qd 8.1117756122 & & 5 & & 10.8952834492 & & 6 & & 12.0523651120 & \cr
& 7 & & 13.4542992617 &
& 8 & & 14.9275108496 & & 9 & & 16.3036898019 & \cr
& 10 & & 17.7350625418 & & 11 & & 18.837088412 & & 12 & & 20.551890978 & \cr
& 13 & & 21.752187480 & & 14 & & 22.93715924&  & 15 & & 23.33859940 & \cr
& 16 & & 23.97767239 & & 17 & & 25.79365179 & & 18 & & 27.1212236 & \cr
& 19 & & 27.8904904 & & 20 & & 28.6462091 & & 21 & & 30.100668 & \cr
& 22 & & 30.884244 & & 23 & & 31.730116 & & 24 & & 32.248613 & \cr
& 25 & & 33.84677 & & 26 & & 34.08053 & & 27 & & 35.12990& \cr
& 28 & & 36.04356 & & 29 & & 36.9637 & & 30 & & 38.2333 & \cr
& 31 & & 39.1512 & & 32 & & 39.7944 & & & & & \cr
\vsp
\noalign{\hrule}
}}$$

\bigskip

$$
\text{Table 2.3}
$$

\redefine\qc{\hskip 20pt}
\halign{$\hfil#\hfil$\quad &$\hfil#\hfil$\quad &$\hfil#$\quad
&$\hfil#\hfil$\quad &$\hfil#$\quad &$\hfil#\hfil$\quad
&$\hfil#$ \cr
\noalign{\hrule}
\noalign{\smallskip}
N & R_0 & I_0 \qc & R_5 & I_5 \qc & R_{10} & I_{10} \qc \cr
\noalign{\smallskip}
\noalign{\hrule}
\noalign{\smallskip}
1000& -3.04& -7.7\times 10^{-2}&
-2.9558& 1.4\times 10^{-2}&
-2.95033& 9.3\times 10^{-4} \cr
2000& -2.85& 3.1\times 10^{-3}&
-2.9629& -3.5\times 10^{-3}&
-2.95651& 5.4\times 10^{-4} \cr
4000& -2.95& -9.7\times 10^{-3}&
-2.9588& 7.3\times 10^{-4}&
-2.95642& -8.7\times 10^{-5} \cr
6000& -2.95& 6.3\times 10^{-2}&
-2.9561& -1.4\times 10^{-3}&
-2.95659& 9.4\times 10^{-5} \cr
8000& -2.98& 6.9\times 10^{-2}&
-2.9558& -9.9\times 10^{-4}&
-2.95660& 4.8\times 10^{-5} \cr
10000& -2.97& 1.8\times 10^{-2}&
-2.9567& -2.3\times 10^{-4}&
-2.95659& 1.1\times 10^{-5} \cr
\noalign{\smallskip}
\noalign{\hrule}
}

\medskip

\redefine\qc{\hskip 23pt}
\redefine\qd{\hskip4.5pt}
\halign{$\hfil#\hfil$\quad &$\hfil#\hfil$\quad &$\hfil#$\quad
&$\hfil#\hfil$\quad &$\hfil#$ \cr
\noalign{\hrule}
\noalign{\smallskip}
N & R_{15} & I_{15} \qc & R_{20} & I_{20} \qc \cr
\noalign{\smallskip}
\noalign{\hrule}
\noalign{\smallskip}
1000& -2.956381& -2.8\times 10^{-3}&
-2.9591449& -4.1\times 10^{-4} \cr
2000& -2.956387& -3.6\times 10^{-4}&
-2.9568292& -9.6\times 10^{-5} \cr
4000& -2.956613& -2.3\times 10^{-5}&
-2.9565854& 4.1\times 10^{-6} \cr
6000& -2.956592& -3.8\times 10^{-6}&
-2.9565951& 2.5\times 10^{-6} \cr
8000& -2.956593& -1.7\times 10^{-6}&
-2.9565940& 6.7\times 10^{-7} \cr
10000& -2.956593& -2.3\times 10^{-6}&
-2.9565932& 1.2\times 10^{-7} \cr
\noalign{\smallskip}
\noalign{\hrule}
}

\medskip

\redefine\qc{\hskip 23pt}
\redefine\qd{\hskip4.5pt}
\halign{$\hfil#\hfil$\quad &$\hfil#\hfil$\quad &$\hfil#$\quad
&$\hfil#\hfil$\quad &$\hfil#$ \cr
\noalign{\hrule}
\noalign{\smallskip}
N & R_{25} & I_{25} \qc & R_{30} & I_{30} \qc \cr
\noalign{\smallskip}
\noalign{\hrule}
\noalign{\smallskip}
1000& -2.9584634& 3.5\times 10^{-3}&
-2.95201141& 3.7\times 10^{-3} \cr
2000& -2.9566395& 1.9\times 10^{-4}&
-2.95638520& 1.7\times 10^{-5} \cr
4000& -2.9565937& -7.5\times 10^{-6}&
-2.95659865& 2.3\times 10^{-6} \cr
6000& -2.9565925& 1.3\times 10^{-7}&
-2.95659358& -2.0\times 10^{-7} \cr
8000& -2.9565932& -1.8\times 10^{-8}&
-2.95659347& 4.3\times 10^{-8} \cr
10000& -2.9565934& -1.0\times 10^{-7}&
-2.95659342& 1.3\times 10^{-8} \cr
\noalign{\smallskip}
\noalign{\hrule}
}

\newpage

\redefine\qc{\hskip 23pt}
\redefine\qd{\hskip4.5pt}
\halign{$\hfil#\hfil$\quad &$\hfil#\hfil$\quad &$\hfil#$ \cr
\noalign{\hrule}
\noalign{\smallskip}
N & R_{35} & I_{35} \qc \cr
\noalign{\smallskip}
\noalign{\hrule}
\noalign{\smallskip}
1000& -2.951829903& 1.3\times 10^{-4} \cr
2000& -2.956587888& -1.9\times 10^{-4} \cr
4000& -2.956589949& 3.4\times 10^{-6} \cr
6000& -2.956593420& 1.6\times 10^{-7} \cr
8000& -2.956593344& 2.0\times 10^{-8} \cr
10000& -2.956593405& -5.0\times 10^{-9} \cr
\noalign{\smallskip}
\noalign{\hrule}
}

\bigskip

$$
\text{Table 2.4}
$$

$$\vbox{\offinterlineskip
\redefine\qa{\hskip6pt}
\redefine\qaa{\hskip7pt}
\redefine\qb{\hskip2.3pt}
\redefine\qc{\hskip16pt}
\redefine\qd{\hskip5.3pt}
\redefine\vsp{height2pt & \omit & & \omit & & \omit & & \omit & & \omit
& & \omit & & \omit & & \omit & \cr}
\define\vspa{height3pt & \omit & & \omit & & \omit & & \omit & & \omit
& & \omit & & \omit & & \omit & \cr}
\halign{&\vrule# &\strut\qa\hfil#\hfil\qb &\vrule#
&\strut\qaa#\hfill\qb \cr
\noalign{\hrule}
\vsp
& n & & \qc $u_n$ & & n & & \qc $u_n$ & & n & & \qc $u_n$ 
& & n & & \qc $u_n$ & \cr
\vsp
\noalign{\hrule}
\vspa
& 1 & & \qd 2.3864500 & & 2 & & \qd 4.3752457 & & 3 & & \qd 6.0435487 & &
4 & & \qd 7.571907 & \cr
& 5 & & \qd 8.841633 & & 6 & & 10.605890 & & 7 & & 11.437474 & &
8 & & 12.76622 & \cr
& 9 & & 13.76869 & & 10 & & 15.2075 & & 11 & & 15.6182 & &
12 & & 16.9663 & \cr
& 13 & & 18.0078 & & 14 & & 18.874 & & & & & & & & & \cr
\vsp
\noalign{\hrule}
}}$$

\bigskip

$$
\text{Table 3.1}
$$

\redefine\qc{\hskip 20pt}
\halign{$\hfil#\hfil$\quad &$\hfil#\hfil$\quad &$\hfil#$\quad
&$\hfil#\hfil$\quad &$\hfil#$\quad &$\hfil#\hfil$\quad
&$\hfil#$ \cr
\noalign{\hrule}
\noalign{\smallskip}
N & R_0 & I_0 \qc & R_5 & I_5 \qc & R_{10} & I_{10} \qc \cr
\noalign{\smallskip}
\noalign{\hrule}
\noalign{\smallskip}
1000& -1.686& 8.7\times 10^{-2}& -1.77582&
1.6 \times 10^{-2}&
-1.782116715& 1.7 \times 10^{-3} \cr
2000& -1.835& -2.9 \times 10^{-2}& -1.78432& -1.4 \times 10^{-3}&
-1.783663963& -1.9 \times 10^{-5} \cr
4000& -1.758& -6.5 \times 10^{-3}& -1.78364& -1.6 \times 10^{-4}&
-1.783642384& -1.0 \times 10^{-6} \cr
6000& -1.820& -5.2 \times 10^{-3}& -1.78366&
2.0 \times 10^{-5}&
-1.783642948& 8.6 \times 10^{-8} \cr
8000& -1.760& -1.1 \times 10^{-2}& -1.78364& -2.0 \times 10^{-5}&
-1.783642816& -2.7 \times 10^{-8} \cr
10000& -1.786& 8.2 \times 10^{-3}& -1.78363& 1.1 \times 10^{-5}&
-1.783642826& 6.5 \times 10^{-9} \cr
\noalign{\smallskip}
\noalign{\hrule}
}

\medskip

\redefine\qc{\hskip 23pt}
\redefine\qd{\hskip4.5pt}
\halign{$\hfil#\hfil$\quad &$\hfil#\hfil$\quad &$\hfil#$\quad
&$\hfil#\hfil$\quad &$\hfil#$ \cr
\noalign{\hrule}
\noalign{\smallskip}
N & R_{15} & I_{15} \qc & R_{20} & I_{20} \qc \cr
\noalign{\smallskip}
\noalign{\hrule}
\noalign{\smallskip}
1000& -1.783149679719& 9.4 \times 10^{-5} \qd & -1.7835293644621857&
-2.4 \times 10^{-5} \qd \cr
2000& -1.783644577360& -3.2 \times 10^{-8} \qd & -1.7836428890128260&
8.1 \times 10^{-9} \qd \cr
4000& -1.783642829247& -2.2 \times 10^{-9} \qd & -1.7836428272010046&
-1.3 \times 10^{-11} \cr
6000& -1.783642827253& 8.0 \times 10^{-11}& -1.7836428271539760&
4.5 \times 10^{-13} \cr
8000& -1.783642827162& -6.1 \times 10^{-12}& -1.7836428271544534&
-7.0 \times 10^{-15} \cr
10000& -1.783642827150& -1.4 \times 10^{-12}& -1.7836428271544125&
-8.3 \times 10^{-16} \cr 
\noalign{\smallskip}
\noalign{\hrule}
}

\newpage

\redefine\qc{\hskip 23pt}
\redefine\qd{\hskip4.5pt}
\halign{$\hfil#\hfil$\quad &$\hfil#\hfil$\quad &$\hfil#$ \cr
\noalign{\hrule}
\noalign{\smallskip}
N & R_{25} & I_{25} \qc \cr
\noalign{\smallskip}
\noalign{\hrule}
\noalign{\smallskip}
1000& -1.7836271941181258672& -1.3 \times 10^{-5} \qd \cr
2000& -1.7836428277732280167& 1.1 \times 10^{-9} \qd \cr
4000& -1.7836428271546296417& 1.4 \times 10^{-15} \cr
6000& -1.7836428271544148349& 9.0 \times 10^{-16} \cr
8000& -1.7836428271544160390& 2.0 \times 10^{-18} \cr
10000& -1.7836428271544160184& -6.2 \times 10^{-19} \cr
\noalign{\smallskip}
\noalign{\hrule}
}

\medskip

\redefine\qc{\hskip 23pt}
\redefine\qd{\hskip4.5pt}
\halign{$\hfil#\hfil$\quad &$\hfil#\hfil$\quad &$\hfil#$ \cr
\noalign{\hrule}
\noalign{\smallskip}
N & R_{30} & I_{30} \qc \cr
\noalign{\smallskip}
\noalign{\hrule}
\noalign{\smallskip}
1000& -1.7836430434865235424107& -3.5 \times 10^{-6} \qd \cr
2000& -1.7836428270926205946250& 5.0 \times 10^{-11} \cr
4000& -1.7836428271544153839298& 1.1 \times 10^{-15} \cr
6000& -1.7836428271544160178686& -1.6 \times 10^{-18} \cr
8000& -1.7836428271544160181358& 2.7 \times 10^{-20} \cr
10000& -1.7836428271544160181708& -3.8 \times 10^{-22} \cr
\noalign{\smallskip}
\noalign{\hrule}
}

\medskip

\redefine\qc{\hskip 23pt}
\redefine\qd{\hskip4.5pt}
\halign{$\hfil#\hfil$\quad &$\hfil#\hfil$\quad &$\hfil#$ \cr
\noalign{\hrule}
\noalign{\smallskip}
N & R_{35} & I_{35} \qc \cr
\noalign{\smallskip}
\noalign{\hrule}
\noalign{\smallskip}
1000& -1.7836435960725853092305493 & -2.3 \times 10^{-7} \qd \cr
2000& -1.7836428271502571300252154 & -3.7 \times 10^{-13} \cr
4000& -1.7836428271544160011102849 & 4.2 \times 10^{-18} \cr
6000& -1.7836428271544160181784305 & -8.5 \times 10^{-21} \cr
8000& -1.7836428271544160181689502 & 1.0 \times 10^{-23} \cr
10000& -1.7836428271544160181690199 & 2.9 \times 10^{-25} \cr
\noalign{\smallskip}
\noalign{\hrule}
}

\bigskip

$$
\text{Table 3.2}
$$

\redefine\qc{\hskip 20pt}
\halign{$\hfil#\hfil$\quad &$\hfil#\hfil$\quad &$\hfil#$\quad
&$\hfil#\hfil$\quad &$\hfil#$\quad &$\hfil#\hfil$\quad
&$\hfil#$ \cr
\noalign{\hrule}
\noalign{\smallskip}
N & R_0 & I_0 \qc & R_5 & I_5 \qc & R_{10} & I_{10} \qc \cr
\noalign{\smallskip}
\noalign{\hrule}
\noalign{\smallskip}
1000& 3.43& -3.9 \times 10^{-1} &
3.3563 & -3.5 \times 10^{-1} &
3.2408605& -3.5 \times 10^{-1} \cr
2000& 3.23& 8.9 \times 10^{-2} &
3.1032& 6.7 \times 10^{-3} &
3.0945490& 2.3 \times 10^{-4} \cr
4000& 3.10& 2.6 \times 10^{-2} &
3.0906& -3.1 \times 10^{-4} &
3.0913811& 2.3 \times 10^{-6} \cr
6000& 3.12& 4.1 \times 10^{-2} &
3.0912& 3.7 \times 10^{-4} &
3.0914333& 2.3 \times 10^{-6} \cr
8000& 3.07& -4.2 \times 10^{-3} &
3.0914& -1.7 \times 10^{-4} &
3.0914340& -3.9 \times 10^{-7} \cr
10000& 3.12& -1.8 \times 10^{-2} &
3.0914& -1.3 \times 10^{-4} &
3.0914342& -1.6 \times 10^{-7} \cr
\noalign{\smallskip}
\noalign{\hrule}
}

\medskip

\redefine\qc{\hskip 23pt}
\redefine\qd{\hskip4.5pt}
\halign{$\hfil#\hfil$\quad &$\hfil#\hfil$\quad &$\hfil#$\quad
&$\hfil#\hfil$\quad &$\hfil#$ \cr
\noalign{\hrule}
\noalign{\smallskip}
N & R_{15} & I_{15} \qc & R_{20} & I_{20} \qc \cr
\noalign{\smallskip}
\noalign{\hrule}
\noalign{\smallskip}
1000& 3.0854245832& -3.1 \times 10^{-1} \qd &
2.923318726655& -2.4 \times 10^{-1} \qd \cr
2000& 3.0919697213& -3.7 \times 10^{-5} \qd &
3.091486059342& -6.5 \times 10^{-5} \qd \cr
4000& 3.0914327625& -3.9 \times 10^{-8} \qd &
3.091434224098& 4.0 \times 10^{-8} \qd \cr
6000& 3.0914342708& 6.1 \times 10^{-8} \qd &
3.091434237097& 5.7 \times 10^{-10} \cr
8000& 3.0914342301& -3.7 \times 10^{-9} \qd &
3.091434236714& -1.0 \times 10^{-11} \cr
10000& 3.0914342356& -3.6 \times 10^{-10} &
3.091434236735& -1.3 \times 10^{-12} \cr
\noalign{\smallskip}
\noalign{\hrule}
}

\medskip

\newpage

\redefine\qc{\hskip 23pt}
\redefine\qd{\hskip4.5pt}
\halign{$\hfil#\hfil$\quad &$\hfil#\hfil$\quad &$\hfil#$\quad
&$\hfil#\hfil$\quad &$\hfil#$ \cr
\noalign{\hrule}
\noalign{\smallskip}
N & R_{25} & I_{25} \qc & R_{30} & I_{30} \qc \cr
\noalign{\smallskip}
\noalign{\hrule}
\noalign{\smallskip}
1000& 2.80762049069122& -1.1 \times 10^{-1} \qd&
2.78668118045678635& 7.4 \times 10^{-2} \qd\cr
2000& 3.09143198592104& -2.2 \times 10^{-5}\qd&
3.09143105396801472& -3.5 \times 10^{-6} \qd\cr
4000& 3.09143423720919& 2.4 \times 10^{-9}\qd&
3.09143423678556380& 3.9 \times 10^{-11} \cr
6000& 3.09143423674457& -5.1 \times 10^{-12}&
3.09143423673870285& -1.7 \times 10^{-13} \cr
8000& 3.09143423673861& 2.3 \times 10^{-13}&
3.09143423673865144& 1.9 \times 10^{-15} \cr
10000& 3.09143423673865& 1.0 \times 10^{-14}&
3.09143423673865098& 6.2 \times 10^{-17} \cr
\noalign{\smallskip}
\noalign{\hrule}
}

\medskip

\redefine\qc{\hskip 23pt}
\redefine\qd{\hskip4.5pt}
\halign{$\hfil#\hfil$\quad &$\hfil#\hfil$\quad &$\hfil#$ \cr
\noalign{\hrule}
\noalign{\smallskip}
N & R_{35} & I_{35} \qc \cr
\noalign{\smallskip}
\noalign{\hrule}
\noalign{\smallskip}
1000& 2.91398635132232829470& 2.7 \times 10^{-1} \qd \cr
2000& 3.09143341359018088974& 1.4 \times 10^{-7} \qd \cr
4000& 3.09143423673979921383& -2.4 \times 10^{-12} \cr
6000& 3.09143423673865010310& -1.1 \times 10^{-15} \cr
8000& 3.09143423673865095777& -1.9 \times 10^{-18} \cr
10000& 3.09143423673865095240& -1.3 \times 10^{-20} \cr
\noalign{\smallskip}
\noalign{\hrule}
}

\bigskip

$$
\text{Table 3.3}
$$

$$\vbox{\offinterlineskip
\redefine\qa{\hskip6pt}
\redefine\qb{\hskip2.3pt}
\redefine\qc{\hskip8pt}
\redefine\qd{\hskip6.4pt}
\redefine\vsp{height2pt & \omit & & \omit & & \omit & & \omit & & \omit
& & \omit & \cr}
\redefine\vspa{height3pt & \omit & & \omit & & \omit & & \omit & & \omit
& & \omit & \cr}
\halign{&\vrule# &\strut\qa\hfil#\hfil\qb \cr
\noalign{\hrule}
\vsp
& n & & \qc $t_n$ & & n & & \qc $t_n$ & & n & & \qc $t_n$ & \cr
\vsp
\noalign{\hrule}
\vspa
& 1 & & \qd 8.2720409199 & & 2 & & 11.3959869930 & & 3 & & 14.8616932015 & \cr
& 4 & & 17.1783243050 & & 5 & & 19.2124566315 & & 6 & & 20.8274294554 & \cr
& 7 & & 23.4659374198 & & 8 & & 25.2726883522 & & 9 & & 27.0035774491 & \cr
& 10 & & 28.1569222690 & & 11 & & 30.2145623343 & & 12 & &
31.6193141164 & \cr
& 13 & & 33.7856279775 & & 14 & & 34.9435854723 & & 15 & &
36.5559515067 & \cr
& 16 & & 37.6356026748 & & 17 & & 39.1608229256 & & 18 & &
40.6589300308 & \cr
& 19 & & 42.8804581030 & & 20 & & 43.2736012304 & & 21 & &
44.9765395474 & \cr
& 22 & & 46.4176568046 & & 23 & & 47.2517710599 & & 24 & &
48.7821808287 & \cr
& 25 & & 50.3519022325 & & 26 & & 51.5688981695 & & 27 & &
53.1356287828 & \cr
& 28 & & 54.0717837181 & & 29 & & 55.0990003336 & & 30 & &
56.4089955139 & \cr
& 31 & & 57.5391214415 & & 32 & & 59.1986375433 & & 33 & &
60.1739007171 & \cr
& 34 & & 61.6441827270 & & 35 & & 62.8146545420 & & 36 & &
63.4247884022 & \cr
& 37 & & 65.1023702197 & & 38 & & 66.0180646898 & & 39 & &
66.8050237006 & \cr
& 40 & & 68.6802278238 & & 41 & & 69.8132058342 & & 42 & &
70.7502185552 & \cr
& 43 & & 71.9861530156 & & 44 & & 72.7927328082 & & 45 & &
74.1137296216 & \cr
& 46 & & 74.8761895173 & & 47 & & 76.2796967025 & & 48 & &
77.4608764665 & \cr
& 49 & & 78.7319975717 & & 50 & & 79.5372511477 & & 51 & &
80.8499015926 & \cr
& 52 & & 81.9286308045 & & 53 & & 82.6995529553 & & 54 & &
83.4681179192 & \cr
& 55 & & 85.2402769759 & & 56 & & 85.6802121224 & & 57 & &
87.2830188249 & \cr
& 58 & & 88.5094955323 & & 59 & & 89.2377130355 & & 60 & &
90.0534073382 & \cr
& 61 & & 91.4472572430 & & 62 & & 92.0496894589 & & 63 & &
93.3566370961 & \cr
& 64 & & 94.2221147468 & & 65 & & 95.3044565474 & & 66 & &
96.6527715250 & \cr
& 67 & & 97.7264314003 & & 68 & & 98.4244540180 & & 69 & &
99.4730638315 & \cr
\vsp
\noalign{\hrule}
}}$$

\newpage

$$
\text{Table 3.4}
$$

$$\vbox{\offinterlineskip
\redefine\qa{\hskip6pt}
\redefine\qb{\hskip2.3pt}
\redefine\qc{\hskip8pt}
\redefine\qd{\hskip6.4pt}
\redefine\vsp{height2pt & \omit & & \omit & & \omit & & \omit & & \omit
& & \omit & \cr}
\redefine\vspa{height3pt & \omit & & \omit & & \omit & & \omit & & \omit
& & \omit & \cr}
\halign{&\vrule# &\strut\qa\hfil#\hfil\qb \cr
\noalign{\hrule}
\vsp
& n & & \qc $u_n$ & & n & & \qc $u_n$ & & n & & \qc $u_n$ & \cr
\vsp
\noalign{\hrule}
\vspa
& 1 & & 12.9399446108 & & 2 & & 15.1248640287 & & 3 & & 17.2775088490 & \cr
& 4 & & 21.9119654118 & & 5 & & 23.7124474310 & & 6 & & 27.6868648494 & \cr
& 7 & & 29.1470584255 & & 8 & & 31.1315265360 & & 9 & & 31.9862854000 & \cr
& 10 & & 33.6323734231 & & 11 & & 35.7361264638 & & 12 & &
38.1008875317 & \cr
& 13 & & 39.9548075690 & & 14 & & 41.3629251312 & & 15 & &
43.0030848131 & \cr
& 16 & & 43.7924301232 & & 17 & & 49.3874980802 & & 18 & &
50.3892911690 & \cr
& 19 & & 51.9883497256 & & 20 & & 53.3610715851 & & 21 & &
55.5058736308 & \cr
& 22 & & 57.1306190068 & & 23 & & 58.5145765119 & & 24 & &
59.2810504632 & \cr
& 25 & & 60.8114113807 & & 26 & & 61.7177742037 & & 27 & &
62.3299217969 & \cr
& 28 & & 65.1200148215 & & 29 & & 66.3871768599 & & 30 & &
67.7658152255 & \cr
& 31 & & 68.5636482576 & & 32 & & 70.0994795387 & & 33 & &
71.9139076205 & \cr
& 34 & & 73.4598285562 & & 35 & & 74.4219698604 & & 36 & &
75.9259426071 & \cr
& 37 & & 76.7219513797 & & 38 & & 80.2370604179 & & 39 & &
80.9795116625 & \cr
& 40 & & 82.2141987387 & & 41 & & 84.1686266809 & & 42 & &
85.4558525934 & \cr
& 43 & & 86.4596989555 & & 44 & & 87.5017176801 & & 45 & &
88.7989284271 & \cr
& 46 & & 90.8930274655 & & 47 & & 91.6094970880 & & 48 & &
93.0648834307 & \cr
& 49 & & 93.8455295242 & & 50 & & 94.4090447654 & & 51 & &
95.8637476000 & \cr
& 52 & & 96.5025102751 & & 53 & & 97.8219970593 & & 54 & &
98.9086789539 & \cr
\vsp
\noalign{\hrule}
}}$$

\bigskip

$$
\text{Table 3.5}
$$

$$\vbox{\offinterlineskip
\redefine\qa{\hskip6pt}
\redefine\qb{\hskip2.3pt}
\redefine\qc{\hskip60pt}
\redefine\qd{\hskip6.4pt}
\redefine\vsp{height2pt & \omit & & \omit & & \omit & & \omit & \cr}
\redefine\vspa{height3pt & \omit & & \omit & & \omit & & \omit & \cr}
\halign{&\vrule# &\strut\qa\hfil#\hfil\qb &\vrule#
&\strut\qa#\hfill\qb \cr
\noalign{\hrule}
\vsp
& n & & \qc $\rho _n$ & & n & & \qc $\rho _n$ & \cr
\vsp
\noalign{\hrule}
\vspa
& 1 & & $3.2308208282+8.9496290911\, i$
& & 2 & & $3.0144204971+19.1670355895\, i$ & \cr
& 3 & & $3.1880664988+26.3033849287\, i$
& & 4 & & $3.1549639910+36.6242398231\, i$ & \cr
& 5 & & $2.7150409653+45.1719799932\, i$
& & 6 & & $2.4938210677+47.5816502442\, i$ & \cr
& 7 & & $3.3624175212+54.4320525502\, i$
& & 8 & & $2.9077749773+64.2513434784\, i$ & \cr
& 9 & & $3.1556119321+71.9344926377\, i$
& & 10 & & $2.4066868777+78.1144688947\, i$ & \cr
& 11 & & $2.9102154501+82.3890698796\, i$
& & 12 & & $2.8870784016+89.7875787849\, i$ & \cr
& 13 & & $3.3672596002+99.9194124003\, i$
& & 14 & & $2.7073645379+107.1592978688\, i$ & \cr
& 15 & & $2.7721770492+110.2613188689\, i$
& & 16 & & $3.1645227049+117.3896482956\, i$ & \cr
& 17 & & $2.5547542302+126.5001914198\, i$
& & 18 & & $2.8669588498+128.2070045099\, i$ & \cr
& 19 & & $3.1444059195+135.3354942155\, i$
& & 20 & & $3.3550056102+145.3938901873\, i$ & \cr
\vsp
\noalign{\hrule}
}}$$

\newpage

$$
\nopagebreak
\vbox{
\hskip 12em \text{Table 3.6}
\vskip 2em
\vbox{\offinterlineskip
\redefine\qa{\hskip6pt}
\redefine\qb{\hskip2.3pt}
\redefine\qc{\hskip8pt}
\redefine\qd{\hskip5.3pt}
\redefine\vsp{height2pt & \omit & & \omit & & \omit & & \omit & & \omit
& & \omit & \cr}
\redefine\vspa{height3pt & \omit & & \omit & & \omit & & \omit & & \omit
& & \omit & \cr}
\halign{&\vrule# &\strut\qa\hfil#\hfil\qb \cr
\noalign{\hrule}
\vsp
& n & & \qc $u_n^+$ & & n & & \qc $u_n^+$ & & n & & \qc $u_n^+$ & \cr
\vsp
\noalign{\hrule}
\vspa
& 1 & & \qd 5.6185671952 & & 2 & & \qd 9.3587692608 & & 3 & & 12.0264936925 & \cr
& 4 & & 13.7962357292 & & 5 & & 16.0894994874 & & 6 & & 17.7847280762 & \cr
& 7 & & 19.4575970308 & & 8 & & 21.8594962784 & & 9 & & 22.5758373316 & \cr
& 10 & & 24.1962387103 & & 11 & & 26.0225517432 & & 12 & &
27.3617234087 & \cr
& 13 & & 29.1281915898 & & 14 & & 29.9334498107 & & 15 & &
31.7848563053 & \cr
& 16 & & 33.0414393853 & & 17 & & 33.9453293541 & & 18 & &
35.7534640159 & \cr
& 19 & & 36.6630191145 & & 20 & & 38.2010379041 & & 21 & &
39.7144586843 & \cr
& 22 & & 40.6779144658 & & 23 & & 41.3507131813 & & 24 & &
43.0237415354 & \cr
& 25 & & 44.5137568744 & & 26 & & 45.2712575072 & & 27 & &
46.4744214908 & \cr
& 28 & & 47.8221666146 & & 29 & & 49.0105055856 & & 30 & &
49.9881789218 & \cr
& 31 & & 51.5155718913 & & 32 & & 51.6631984921 & & 33 & &
53.5375246718 & \cr
& 34 & & 54.4042145877 & & 35 & & 55.5750957955 & & 36 & &
58.8270009777 & \cr
& 37 & & 59.7418610167 & & 38 & & 61.2006909669 & & 39 & &
61.7440856002 & \cr
& 40 & & 63.2789136582 & & 41 & & 64.1010371646 & & 42 & &
64.8912913186 & \cr
& 43 & & 66.3286609040 & & 44 & & 67.5430514279 & & 45 & &
70.1770012608 & \cr
& 46 & & 71.4492522675 & & 47 & & 72.3027663172 & & 48 & &
73.1469069994 & \cr
& 49 & & 74.5424808567 & & 50 & & 75.2004527728 & & 51 & &
76.3094730233 & \cr
& 52 & & 77.0682809138 & & 53 & & 78.2931224271 & & 54 & &
79.2581319133 & \cr
& 55 & & 80.0330330946 & & 56 & & 81.1315223512 & & 57 & &
81.9704731391 & \cr
& 58 & & 83.0719604785 & & 59 & & 84.1964169499 & & 60 & &
85.2011952740 & \cr
& 61 & & 85.7778751490 & & 62 & & 86.9407183317 & & 63 & &
87.4479355411 & \cr
& 64 & & 88.7059787721 & & 65 & & 89.9368381341 & & 66 & &
90.6181312264 & \cr
& 67 & & 91.4583733416 & & 68 & & 92.7393858316 & & 69 & &
93.0881396596 & \cr
& 70 & & 94.5091774704 & & 71 & & 95.4769845208 & & 72 & &
96.1116129723 & \cr
& 73 & & 96.9317027531 & & 74 & & 97.8372552255 & & 75 & &
99.1010053953 & \cr
& 76 & & 99.8881597950 & & & & & & & & & \cr
\vsp
\noalign{\hrule}
\vsp
& n & & \qc $u_n^-$ & & n & & \qc $u_n^-$ & & n & & \qc $u_n^-$ & \cr
\vsp
\noalign{\hrule}
\vspa
& 1 & &\hskip -48pt 0 & & 2 & & 24.4022873037 & & 3 & & 26.6418851276 & \cr
& 4 & & 29.3670678246 & & 5 & & 33.5747954436 & & 6 & & 35.2863556538 & \cr
& 7 & & 38.5418813017 & & 8 & & 40.0447318001 & & 9 & & 44.8539372315 & \cr
& 10 & & 46.8484465576 & & 11 & & 50.0699839799 & & 12 & &
52.1256566323 & \cr
& 13 & & 54.2746979473 & & 14 & & 55.4033486176 & & 15 & &
56.5807632829 & \cr
& 16 & & 58.8640673277 & & 17 & & 60.9985910184 & & 18 & &
63.8456398387 & \cr
& 19 & & 65.1978869599 & & 20 & & 71.4691922483 & & 21 & &
72.8978373808 & \cr
& 22 & & 75.0463737748 & & 23 & & 76.6579317410 & & 24 & &
78.6419254744 & \cr
& 25 & & 80.0462305996 & & 26 & & 83.8815457436 & & 27 & &
85.2027989008 & \cr
& 28 & & 86.6795992346 & & 29 & & 88.7015447955 & & 30 & &
90.5599843766 & \cr
& 31 & & 92.9081333908 & & 32 & & 94.5315054431 & & 33 & &
96.4777165784 & \cr
& 34 & & 97.7100184923 & & 35 & & 99.3454122963 & & & & & \cr
\vsp
\noalign{\hrule}
}}}
$$

\newpage

$$
\text{Table 3.7}
$$

$$\vbox{\offinterlineskip
\redefine\qa{\hskip6pt}
\redefine\qb{\hskip2.3pt}
\redefine\qc{\hskip60pt}
\redefine\qd{\hskip6.4pt}
\redefine\vsp{height2pt & \omit & & \omit & & \omit & & \omit & \cr}
\redefine\vspa{height3pt & \omit & & \omit & & \omit & & \omit & \cr}
\halign{&\vrule# &\strut\qa\hfil#\hfil\qb &\vrule#
&\strut\qa#\hfill\qb \cr
\noalign{\hrule}
\vsp
& n & & \qc $\rho _n^+$ & & n & & \qc $\rho _n^+$ & \cr
\vsp
\noalign{\hrule}
\vspa
& 1 & & $3.3591319232+57.3250633340\, i$
& & 2 & & $3.5001503209+68.5679322965\, i$ & \cr
\vsp
\noalign{\hrule}
\vsp
& n & & \qc $\rho _n^-$ & & n & & \qc $\rho _n^-$ & \cr
\vsp
\noalign{\hrule}
\vspa
& 1 & & $5.7692647648+8.9956889852\, i$
& & 2 & & $4.8476735625+14.0858508094\, i$ & \cr
& 3 & & $5.3846794177+18.2757545274\, i$
& & 4 & & $4.2067408135+20.6821222248\, i$ & \cr
& 5 & & $5.3861276711+26.6587658619\, i$
& & 6 & & $4.5402410961+31.5997171480\, i$ & \cr
& 7 & & $5.7862146670+36.6524415925\, i$
& & 8 & & $4.3281324452+41.9810464020\, i$ & \cr
& 9 & & $5.2171816339+44.9332286844\, i$
& & 10 & & $4.7664043254+48.8240853600\, i$ & \cr
& 11 & & $5.6527005509+54.2712676115\, i$
& & 12 & & $4.6702375442+60.1227186672\, i$ & \cr
& 13 & & $5.2007749392+63.8919812818\, i$
& & 14 & & $3.7792244502+66.9023711945\, i$ & \cr
& 15 & & $3.5790185110+69.3185823254\, i$
& & 16 & & $5.7686816361+72.2279984677\, i$ & \cr
& 17 & & $4.6599017691+76.8392451729\, i$
& & 18 & & $5.3421407237+81.9810448650\, i$ & \cr
& 19 & & $3.3026675429+82.2734272886\, i$
& & 20 & & $3.9986981918+88.8368300771\, i$ & \cr
& 21 & & $4.9626158230+90.6464501774\, i$
& & 22 & & $4.7401362618+94.8547536615\, i$ & \cr
& 23 & & $5.8412001878+99.7795574645\, i$ & & & & & \cr
\vsp
\noalign{\hrule}
}}$$

\bigskip

$$
\text{Table 4.1a}
$$

\redefine\qc{\hskip 20pt}
\halign{$\hfil#\hfil$\quad &$\hfil#\hfil$\quad &$\hfil#$\quad
&$\hfil#\hfil$\quad &$\hfil#$\quad &$\hfil#\hfil$\quad
&$\hfil#$ \cr
\noalign{\hrule}
\noalign{\smallskip}
N & R_0 & I_0 \qc & R_1 & I_1 \qc & R_2 & I_2 \qc \cr
\noalign{\smallskip}
\noalign{\hrule}
\noalign{\smallskip}
1000& 2.22& 6.3\times 10^{-3}&
2.1777& -1.1\times 10^{-2}&
2.173977& 2.6\times 10^{-4} \cr
5000& 2.14& -2.2\times 10^{-2}&
2.1700& 3.7\times 10^{-3}&
2.173544& -1.2\times 10^{-4} \cr
10000& 2.18& -1.3 \times 10^{-2}&
2.1753& 1.1 \times 10^{-3}&
2.173833& -1.0\times 10^{-4} \cr
30000& 2.19& -2.2\times 10^{-2}&
2.1729& -6.9\times 10^{-4}&
2.173759& -1.0\times 10^{-5} \cr
100000& 2.17& 1.2\times 10^{-2}&
2.1735& 1.0\times 10^{-4}&
2.173747& -4.7\times 10^{-6} \cr
\noalign{\smallskip}
\noalign{\hrule}
}

\medskip

\redefine\qc{\hskip 20pt}
\halign{$\hfil#\hfil$\quad &$\hfil#\hfil$\quad &$\hfil#$\quad
&$\hfil#\hfil$\quad &$\hfil#$\quad &$\hfil#\hfil$\quad
&$\hfil#$ \cr
\noalign{\hrule}
\noalign{\smallskip}
N & R_3 & I_3 \qc & R_4 & I_4 \qc & R_5 & I_5 \qc \cr
\noalign{\smallskip}
\noalign{\hrule}
\noalign{\smallskip}
1000& 2.167862& -6.2\times 10^{-3}&
2.20484& -5.4\times 10^{-2}&
2.49186& 7.1\times 10^{-2} \cr
5000& 2.173880& -6.8\times 10^{-4}&
2.17935& -1.6\times 10^{-4}&
2.18302& 2.7\times 10^{-2} \cr
10000& 2.173937& 1.8\times 10^{-4}&
2.17267& 1.6\times 10^{-3}&
2.16380& -2.8\times 10^{-3} \cr
30000& 2.173700& -2.3\times 10^{-6}&
2.17367& -3.7\times 10^{-4}&
2.17551& -9.3\times 10^{-4} \cr
100000& 2.173745& -5.9\times 10^{-6}&
2.17378& -5.1\times 10^{-5}&
2.17406& 9.7\times 10^{-5} \cr
\noalign{\smallskip}
\noalign{\hrule}
}

\newpage

\redefine\qc{\hskip 20pt}
\redefine\qd{\hskip 21pt}
\redefine\qe{\hskip 7.5pt}
\halign{$\hfil#\hfil$\quad &$\hfil#\hfil$\quad &$\hfil#$\quad
&$\hfil#\hfil$\quad &$\hfil#$\quad &$\hfil#\hfil$\quad
&$\hfil#$ \cr
\noalign{\hrule}
\noalign{\smallskip}
N & R_6 & I_6 \qc & R_7 & I_7 \qc & R_8 & I_8 \qc \cr
\noalign{\smallskip}
\noalign{\hrule}
\noalign{\smallskip}
1000& 2.352& 1.2 \qd &
-1.442 \qe & 2.2 \qd &
-7.975 \qe & -6.1 \qd \cr
5000& 2.082& 7.3\times 10^{-2}&
1.835& -1.7\times 10^{-1}&
2.195& -1.0 \qd \cr
10000& 2.170& -4.1\times 10^{-2}&
2.293& -6.2\times 10^{-2}&
2.480& 2.2\times 10^{-1} \cr
30000& 2.179& 5.3\times 10^{-3}&
2.165& 2.4\times 10^{-2}&
2.101& 1.0\times 10^{-2} \cr
100000& 2.173& 1.3\times 10^{-3}&
2.169& 1.8\times 10^{-3}&
2.164& -7.2\times 10^{-3} \cr
\noalign{\smallskip}
\noalign{\hrule}
}

\bigskip

$$
\text{Table 4.1b}
$$

\redefine\qc{\hskip 20pt}
\halign{$\hfil#\hfil$\quad &$\hfil#\hfil$\quad &$\hfil#$\quad
&$\hfil#\hfil$\quad &$\hfil#$\quad &$\hfil#\hfil$\quad
&$\hfil#$ \cr
\noalign{\hrule}
\noalign{\smallskip}
N & R_0 & I_0 \qc & R_1 & I_1 \qc & R_2 & I_2 \qc \cr
\noalign{\smallskip}
\noalign{\hrule}
\noalign{\smallskip}
1000& 3.263& 8.3\times 10^{-2}&
3.3048& 6.1\times 10^{-2}&
3.28176& 2.6\times 10^{-2} \cr
5000& 3.305& -1.1\times 10^{-2}&
3.2708& -1.8\times 10^{-2}&
3.26879& -8.6\times 10^{-4} \cr
10000& 3.245& 6.4\times 10^{-3}&
3.2612& -4.8\times 10^{-3}&
3.26487& 6.8\times 10^{-4} \cr
30000& 3.263& 2.9\times 10^{-2}&
3.2697& 9.0\times 10^{-4}&
3.26616& 1.7\times 10^{-4} \cr
100000& 3.276& -8.2\times 10^{-3}&
3.2665& 5.8\times 10^{-4}&
3.26610& 4.9\times 10^{-5} \cr
\noalign{\smallskip}
\noalign{\hrule}
}

\medskip

\redefine\qc{\hskip 20pt}
\redefine\qd{\hskip 21pt}
\redefine\qe{\hskip 7.5pt}
\halign{$\hfil#\hfil$\quad &$\hfil#\hfil$\quad &$\hfil#$\quad
&$\hfil#\hfil$\quad &$\hfil#$\quad &$\hfil#\hfil$\quad
&$\hfil#$ \cr
\noalign{\hrule}
\noalign{\smallskip}
N & R_3 & I_3 \qc & R_4 & I_4 \qc & R_5 & I_5 \qc \cr
\noalign{\smallskip}
\noalign{\hrule}
\noalign{\smallskip}
1000& 3.38022& 2.5\times 10^{-2}&
3.524& 2.2 \qd &
-33.06 \qe & 8.3 \qd \cr
5000& 3.26878& 6.2\times 10^{-3}&
3.080& 7.2\times 10^{-2}&
1.75& -2.9 \qd \cr
10000& 3.26486& -2.7\times 10^{-3}&
3.329& -3.1\times 10^{-2}&
3.88& 1.0 \qd \cr
30000& 3.26672& -1.6\times 10^{-4}&
3.271& 1.2 \times 10^{-2}&
3.06& 1.0\times 10^{-1} \cr
100000& 3.26609& 6.8\times 10^{-5}&
3.264& -1.3\times 10^{-3}&
3.28& -3.1\times 10^{-2} \cr
\noalign{\smallskip}
\noalign{\hrule}
}

\bigskip

$$
\text{Table 4.2}
$$

$$\vbox{\offinterlineskip
\redefine\qa{\hskip6pt}
\redefine\qaa{\hskip7pt}
\redefine\qb{\hskip2.3pt}
\redefine\qc{\hskip16pt}
\redefine\qd{\hskip5pt}
\redefine\vsp{height2pt & \omit & & \omit & & \omit & & \omit & & \omit
& & \omit & & \omit & & \omit & \cr}
\redefine\vspa{height3pt & \omit & & \omit & & \omit & & \omit & & \omit
& & \omit & & \omit & & \omit & \cr}
\halign{&\vrule# &\strut\qa$\hfil#\hfil$\qb &\vrule#
&\strut\qaa$#\hfill$\qb \cr
\noalign{\hrule}
\vsp
& n & & \qc t_n & & n & & \qc t_n & & n & & \qc t_n 
& & n & & \qc t_n & \cr
\vsp
\noalign{\hrule}
\vspa
& 1 & & 10.2562 & & 2 & & 13.6866 & & 3 & & 15.9599 & &
4 & & 17.038 & \cr
& 5 & & 19.026 & & 6 & & 20.017 & & 7 & & 22.472 & &
8 & & 23.745 & \cr
& 9 & & 25.351 & & 10 & & 26.229 & & 11 & & 27.561 & &
12 & & 28.847 & \cr
& 13 & & 29.986 & & & & & & & & & & & & & \cr
\vsp
\noalign{\hrule}
\vspa
& -1 & & \qd -3.12740 & & -2 & & \qd  -6.5577 & & -3 & & \qd -8.8310 & &
-4 & &  \qd -9.9092 & \cr
& -5 & & -11.8976 & & -6 & & -12.888 & & -7 & & -15.343 & &
-8 & & -16.616 & \cr
& -9 & & -18.222 & & -10 & & -19.101 & & -11 & & -20.433 & &
-12 & & -21.718 & \cr
& -13 & & -22.857 & & -14 & & -24.859 & & -15 & & -25.892 & &
-16 & & -26.850 & \cr
& -17 & & -28.418 & & -18 & & -28.927 & & & & & & & & & \cr
\vsp
\noalign{\hrule}
}}$$

\newpage

$$
\text{Table 4.3}
$$

$$\vbox{\offinterlineskip
\redefine\qa{\hskip6pt}
\redefine\qb{\hskip2.3pt}
\redefine\qc{\hskip8pt}
\redefine\qd{\hskip6.4pt}
\redefine\vsp{height2pt & \omit & & \omit & & \omit & & \omit & & \omit
& & \omit & & \omit & & \omit & \cr}
\redefine\vspa{height3pt & \omit & & \omit & & \omit & & \omit & & \omit
& & \omit & & \omit & & \omit & \cr}
\halign{&\vrule# &\strut\qa$\hfil#\hfil$\qb \cr
\noalign{\hrule}
\vsp
& k & & n & & T & & N(T) & & k & & n & & T & & N(T) & \cr
\vsp
\noalign{\hrule}
\vspa
& \bold Q (\sqrt {2}) & & 1 & & 100 & & 84
& & \bold Q(\sqrt {2}) & & -1 & & 100 & & 93 & \cr
& \bold Q (\sqrt {2}) & & 2 & & 100 & & 82
& & \bold Q (\sqrt {2}) & & -2 & & 100 & & 96 & \cr
& \bold Q (\sqrt {2}) & & 3 & & 100 & & 81
& & \bold Q (\sqrt {2}) & & -3 & & 100 & & 98 & \cr
& \bold Q (\sqrt {2}) & & 4 & & 100 & & 80
& & \bold Q (\sqrt {2}) & & -4 & & 100 & & 100 & \cr
& \bold Q (\sqrt {2}) & & 5 & & 100 & & 78
& & \bold Q (\sqrt {2}) & & -5 & & 100 & & 101 & \cr
& \bold Q (\sqrt {2}) & & 10 & & 100 & & 79
& & \bold Q (\sqrt {2}) & & -10 & & 100 & & 108 & \cr
& \bold Q (\sqrt {5}) & & 1 & & 100 & & 75
& & \bold Q (\sqrt {5}) & & -1 & & 100 & & 88 & \cr
& \bold Q (\sqrt {5}) & & 2 & & 100 & & 72
& & \bold Q (\sqrt {5}) & & -2 & & 100 & & 92 & \cr
& \bold Q (\sqrt {5}) & & 3 & & 100 & & 71
& & \bold Q (\sqrt {5}) & & -3 & & 100 & & 95 & \cr
& \bold Q (\sqrt {5}) & & 4 & & 100 & & 71
& & \bold Q(\sqrt {5}) & & -4 & & 100 & & 97 & \cr
& \bold Q (\sqrt {5}) & & 5 & & 100 & & 71
& & \bold Q(\sqrt {5}) & & -5 & & 100 & & 99 & \cr
& \bold Q (\sqrt {5}) & & 10 & & 100 & & 93
& & \bold Q(\sqrt {5}) & & -10 & & 100 & & 106 & \cr
& \bold Q (\sqrt {19}) & & 1 & & 50 & & 51
& & \bold Q (\sqrt {19}) & & -1 & & 50 & & 52 & \cr
& \bold Q (\sqrt {19}) & & 2 & & 50 & & 50
& & \bold Q (\sqrt {19}) & & -2 & & 50 & & 53 & \cr
& \bold Q (\sqrt {29}) & & 1 & & 100 & & 107
& & \bold Q (\sqrt {29}) & & -1 & & 100 & & 112 & \cr
& \bold Q (\sqrt {29}) & & 2 & & 100 & & 105
& & \bold Q (\sqrt {29}) & & -2 & & 100 & & 114 & \cr
& \bold Q (\sqrt {29}) & & 3 & & 100 & & 104
& & \bold Q (\sqrt {29}) & & -3 & & 100 & & 116 & \cr
& \bold Q (\sqrt {29}) & & 4 & & 100 & & 102
& & \bold Q (\sqrt {29}) & & -4 & & 100 & & 117 & \cr
& \bold Q (\sqrt {29}) & & 5 & & 100 & & 102
& & \bold Q (\sqrt {29}) & & -5 & & 100 & & 118 & \cr
& \bold Q (\sqrt {29}) & & 10 & & 100 & & 99
& & \bold Q (\sqrt {29}) & & -10 & & 100 & & 123 & \cr
& \bold Q (\sqrt {31}) & & 1 & & 40 & & 41
& & \bold Q (\sqrt {31}) & & -1 & & 40 & & 42 & \cr
& \bold Q (\sqrt {67}) & & 1 & & 30 & & 32
& & \bold Q (\sqrt {67}) & & -1 & & 30 & & 33 & \cr
\vsp
\noalign{\hrule}
}}$$

\newpage

$$
\text{Table 5.1}
$$

\redefine\qc{\hskip 20pt}
\halign{$\hfil#\hfil$\quad &$\hfil#\hfil$\quad &$\hfil#$\quad
&$\hfil#\hfil$\quad &$\hfil#$\quad &$\hfil#\hfil$\quad
&$\hfil#$ \cr
\noalign{\hrule}
\noalign{\smallskip}
N & R_0 & I_0 \qc & R_2 & I_2 \qc & R_4 & I_4 \qc \cr
\noalign{\smallskip}
\noalign{\hrule}
\noalign{\smallskip}
1000& 3.09& 5.3\times 10^{-1}&
3.334& 4.5\times 10^{-1}&
3.483& 4.9\times 10^{-1} \cr
5000& 3.31& -5.5\times 10^{-2}&
3.317& -2.4\times 10^{-2}&
3.335& -3.2\times 10^{-3} \cr
10000& 3.42& 1.5 \times 10^{-1}&
3.417& -7.8 \times 10^{-3}&
3.389& -1.3\times 10^{-2} \cr
30000& 3.45& -1.4\times 10^{-1}&
3.385& -2.3\times 10^{-2}&
3.380& -6.2\times 10^{-3} \cr
100000& 3.25& -4.8\times 10^{-2}&
3.375& 1.6\times 10^{-3}&
3.383& 1.2\times 10^{-3} \cr
\noalign{\smallskip}
\noalign{\hrule}
}

\medskip

\redefine\qc{\hskip 20pt}
\halign{$\hfil#\hfil$\quad &$\hfil#\hfil$\quad &$\hfil#$\quad
&$\hfil#\hfil$\quad &$\hfil#$\quad &$\hfil#\hfil$\quad
&$\hfil#$ \cr
\noalign{\hrule}
\noalign{\smallskip}
N & R_6& I_6 \qc & R_8 & I_8 \qc & R_{10} & I_{10} \qc \cr
\noalign{\smallskip}
\noalign{\hrule}
\noalign{\smallskip}
1000& 3.9634& 6.4\times 10^{-1}&
4.03999& 6.2\times 10^{-1}&
4.792707& 5.2\times 10^{-1} \cr
5000& 3.3419& 3.6\times 10^{-2}&
3.38817& 4.1\times 10^{-2}&
3.400866& 6.3\times 10^{-2} \cr
10000& 3.3825& -1.0\times 10^{-2}&
3.37107& -4.5\times 10^{-3}&
3.377550& 4.8\times 10^{-3} \cr
30000& 3.3803& 1.9\times 10^{-4}&
3.38367& 1.4\times 10^{-3}&
3.383935& 3.6\times 10^{-4} \cr
100000& 3.3839& 1.2\times 10^{-4}&
3.38369& -8.5\times 10^{-5}&
3.383657& -8.2\times 10^{-6} \cr
\noalign{\smallskip}
\noalign{\hrule}
}

\bigskip

$$
\text{Table 5.2}
$$

$$\vbox{\offinterlineskip
\redefine\qa{\hskip6pt}
\redefine\qaa{\hskip7pt}
\redefine\qb{\hskip2.3pt}
\redefine\qc{\hskip10pt}
\define\qca{\hskip14pt}
\define\qcb{\hskip12pt}
\define\qcc{\hskip8pt}
\redefine\qd{\hskip5.3pt}
\redefine\vsp{height2pt & \omit & & \omit & & \omit & & \omit & & \omit
& & \omit & & \omit & & \omit & & \omit & & \omit & \cr}
\define\vspa{height3pt & \omit & & \omit & & \omit & & \omit & & \omit
& & \omit & & \omit & & \omit & & \omit & & \omit & \cr}
\halign{&\vrule# &\strut\qa\hfil#\hfil\qb &\vrule#
&\strut\qaa#\hfill\qb \cr
\noalign{\hrule}
\vsp
& n & & \qca $u_n$ & & n & & \qcb $u_n$ & & n & & \qc $u_n$ 
& & n & & \qc $u_n$ & & n & & \qcc $u_n$ & \cr
\vsp
\noalign{\hrule}
\vspa
& 1 & & \qd 2.79373 & & 2 & & \qd 4.0887 & & 3 & & \qd 5.362 & &
4 & & \qd 5.887 & & 5 & & \qd 7.03 & \cr
& 6 & & \qd 7.46 & & 7 & & \qd 7.90 & &
8 & & \qd 8.8 & & 9 & & \qd 9.6 & & & & & \cr
\vsp
\noalign{\hrule}
}}$$

\bigskip

$$
\text{Table 6.1}
$$

$$\vbox{\offinterlineskip
\redefine\qa{\hskip6pt}
\redefine\qaa{\hskip7pt}
\redefine\qb{\hskip2.3pt}
\redefine\qc{\hskip10pt}
\redefine\qca{\hskip14pt}
\redefine\qcb{\hskip12pt}
\redefine\qcc{\hskip8pt}
\redefine\qd{\hskip5.3pt}
\redefine\vsp{height2pt & \omit & & \omit & & \omit & & \omit & & \omit
& & \omit & & \omit & & \omit & & \omit & \cr}
\define\vspa{height3pt & \omit & & \omit & & \omit & & \omit & & \omit
& & \omit & & \omit & & \omit & & \omit & \cr}
\halign{&\vrule# &\strut\qa\hfil#\hfil\qb &\vrule# 
&\strut\qa\hfill#\hfill\qb \cr
\noalign{\hrule}
\vsp
& N & & $r_N^{(1)}$ & & $r_N^{(2)}$ & & $r_N^{(3)}$ 
& & $r_N^{(4)}$ & & N
& & $r_N^{(5)}$ & & N & & $r_N^{(6)}$ & \cr
\vsp
\noalign{\hrule}
\vspa
& 5 & & 0.96 & & 0.98 & & 0.93 & & 0.92 & & 0 & &
1.03 & & 0 & & 0.98 & \cr
& 10 & & 0.95 & & 1.05 & & 1.12 & & 0.84 & & 2 & &
0.82 & & 1 & & 0.99 & \cr
& 15 & & 1.07 & & 1.06 & & 0.85 & & 1.03 & & 4 & &
0.81 & & 2 & & 1.01 & \cr
& 20 & & 1.07 & & 1.09 & & 1.23 & & 1.07 & & 6 & &
1.19 & & 3 & & 0.99 & \cr
& 25 & & 1.06 & & 0.90 & & 1.25 & & 0.87 & & 8 & &
1.39 & & 4 & & 1.01 & \cr
& 30 & & 0.88 & & 0.88 & & 0.80 & & 1.22 & & 10 & &
0.67 & & 5 & & 0.94 & \cr
\vsp
\noalign{\hrule}
}}$$

\input epsf.tex
\midinsert
\topcaption{ Figure 7.1 }
\endcaption
%\vskip 3em
\centerline{
\epsfxsize=1.0\hsize
\epsfysize=1.0\vsize
\epsfbox{g2.ps}
}
\endinsert

\newpage
\pageno=36

\centerline{References}

\bigskip

\item{[E]} H.M. Edwards, Riemann's zeta function, Academic Press, 1974.

\medskip

\item {[In]} A.E. Ingham, The distribution of prime numbers,
Cambridge mathematical library series, 1990.

\medskip

\item {[Is]} H. Ishii, On calculations of zeros of $L$-functions associated
with cusp forms, Memoirs Inst. Science and Engineering, 
Ritsumeikan Univ. 50(1991), 163--172 (in Japanese).

\medskip

\item {[H]} E. Hecke, \"Uber analytische Fuktionen und die Verteilung
von Zahlen mod. eins, Hamburg Abhandlungen, 1(1921), 54--76
(=Werke No. 16).

\medskip

\item{[KZ]} W. Kohnen and D. Zagier, Values of $L$-series
of modular forms at the center of the critical strip, Inv. Math.
64(1981), 175--198.

\medskip

\item{[L]} S. Lang, Algebraic number theory, Addison-Wesley, 1970.

\medskip

\item{[LL]} R.P. Langlands, On the functional equation of Artin
$L$-functions, Yale University Lecture note.

\medskip

\item {[M]} J.S. Milne, Motives over finite fields, Proc. Symposia Pure 
Math. 55(1994), part 1, 401--459.

\medskip

\item{[Se]} J-P. Serre, Une interpr\'etation des congruences
relatives \`a la fonctions $\tau$ de Ramanujan, S\'eminaires
Delange-Pisot-Poitou 1967/68, $\text {n}^{\text {o}}$ 14.

\medskip

\item {[Sha1]} F. Shahidi, Third symmetric power $L$-functions for $GL(2)$,
Comp. Math. 70
\newline (1989), 245--273.

\medskip

\item {[Sha2]} F. Shahidi, Symmetric power $L$-functions for $GL(2)$,
CRM Proceedings \& Lecture notes 4(1994), 159--182.

\medskip

\item{[Sh1]} G. Shimura, Introduction to the arithmetic theory of
automorphic functions, Iwanami-Shoten and Princeton University Press, 1971.

\medskip

\item{[Sh2]} G. Shimura, On modular forms of half integral weight,
Ann. of Math. 97(1973), 440--481.

\medskip

\item{[W]} A. Weil, Basic number theory, Die Grundlehren der
mathematischen Wissenscha-
\newline ften 144, Springer Verlag, 1967.

\medskip

\item{[WW]} E.T. Whittaker and G.N. Watson, A course of Modern analysis,
fourth edition, Cambridge University Press, 1927.

\medskip

\item{[Y1]} H. Yoshida, On a certain distribution on $GL(n)$ and explicit
formulas, Proc. Japan Acad. Ser. A, 63(1987), 396--399.

\medskip

\item{[Y2]} H. Yoshida, On calculations of zeros of $L$-functions
related with Ramanujan's discriminant function on the critical line,
J. of Ramanujan Math. Soc., Ramanujan Birth Centenary Special
Issue, 3(1988), 87--95.

\medskip

\item{[Y3]} H. Yoshida, On hermitian forms attached to zeta functions,
Adv. Stud. in pure math. 21(1992), 281--325.

\enddocument